\documentclass{amsart}
\usepackage{amsmath,amssymb,amsfonts}
\usepackage{colordvi}
\usepackage[all]{xy}

\def\mpn{\medskip\par\noindent}

\def\mmpn{\vskip 1em minus 1em\par\noindent}

\def\smp{\smallskip\par}

\def\CB{{\mathcal B}}
\def\CE{{\mathcal E}}
\def\CF{{\mathcal F}}
\def\CO{{\mathcal O}}
\def\CP{{\mathcal P}}
\def\CR{{\mathcal R}}
\def\CS{{\mathcal S}}
\def\CW{{\mathcal W}}

\def\Ker{\operatorname{Ker}\nolimits}

\def\Id{\operatorname{id}\nolimits}
\def\Res{\operatorname{Res}\nolimits}
\def\Ind{\operatorname{Ind}\nolimits}

\def\Card{\operatorname{Card}\nolimits}
\def\tr{\operatorname{tr}\nolimits}

\def\ls#1#2{{\,^{#1}\!#2}}

\def\Z{\mathbb{Z}}
\def\N{\mathbb{N}}
\def\Q{\mathbb{Q}}

\def\pf{\par\bigskip\noindent{\bf Proof~: }}
\def\endpf{\nolinebreak~\leaders\hbox to 1em{\hss\
\hss}\hfill~\raisebox{.5ex}
{\framebox[1ex]{}}\par\bigskip}

\makeatletter
\renewenvironment{enumerate}{\ifnum \@enumdepth >3 \@toodeep\else
       \advance\@enumdepth \@ne
       \edef\@enumctr{enum\romannumeral\the\@enumdepth}\list
       {\csname  label\@enumctr\endcsname}{\setlength{\topsep}{1ex}
\setlength{\itemsep}{0 pt}\usecounter
         {\@enumctr}\def\makelabel##1{\hss\llap{##1}}}\fi}{\endlist}

\renewenvironment{itemize}{\ifnum \@itemdepth >3 \@toodeep\else
\advance\@itemdepth \@ne
\edef\@itemitem{labelitem\romannumeral\the\@itemdepth}
\list{\csname\@itemitem\endcsname}{\setlength{\topsep}{1ex}\setlength
{\itemsep}{0pt}\def\makelabel##1{\hss\llap{##1}}}\fi}
{\endlist}
\def\@seccntformat#1{\csname the#1\endcsname.\quad}

\def\section{\pagebreak[3]\setcounter{prop}{0}\setcounter{equation}{0}
\@startsection{section}{1}{\z@}{4ex plus  6ex}{2ex}{\center\reset@font \large\bf}}

\makeatother

\def\theprop{\thesection.\arabic{prop}}

\newenvironment{enonce}[1]{\pagebreak[3]\refstepcounter{prop}\mmpn
{{\bf  \thesection.\arabic{prop}.\ #1.}}\begin{it} }{\end{it}\smp}
\def\thesection{\arabic{section}}

\newcommand{\result}[1]{\begin{enonce}{#1}}
\newcommand{\fresult}{\end{enonce}}

\begin{document}

\title[The algebra of essential relations on a finite set]
{The algebra of essential relations on a finite set} 

\author{Serge Bouc}
\author{Jacques Th\'evenaz}
\date\today


\begin{abstract}
Let $X$ be a finite set and let $k$ be a commutative ring.
We consider the $k$-algebra of the monoid of all relations on~$X$,
modulo the ideal generated by the relations factorizing through a set of cardinality strictly smaller than $\Card(X)$,
called inessential relations. This quotient is called the essential algebra associated to $X$.
We then define a suitable nilpotent ideal of the essential algebra
and describe completely the structure of the corresponding quotient,
a product of matrix algebras over suitable group algebras.
In particular, we obtain a description of the Jacobson radical
and of all the simple modules for the essential algebra.
\end{abstract}

\maketitle


\section{Introduction}

\noindent
Let $X$ and $Y$ be finite sets. A {\it correspondence\/} between $X$ and~$Y$ is a subset $R$ of $X\times Y$.
In case $X=Y$, we say that $R$ is a {\it relation\/} on~$X$.
Correspondences can be composed as follows. If $R\subseteq X\times Y$ and $S\subseteq Y\times Z$,
then $RS$ is the correspondence between $X$ and~$Z$ defined by
$$RS=\{ (x,z)\in X\times Z \,\mid\, \exists \, y\in Y \;\text{ such that } \; (x,y)\in R  \;\text{ and } \; (y,z)\in S \} \,.$$
In particular the set of all relations on~$X$ is a monoid.
Given a commutative ring $k$ and a finite set $X$,
let $\CR$ be the $k$-algebra of the monoid of all relations on~$X$ (having this monoid as a $k$-basis).

Throughout this paper, $X$ will denote a finite set of cardinality~$n$.
We say that a relation $R$ on~$X$ is {\it inessential\/} if there exists a set $Y$ with $\Card(Y)<\Card(X)$
and two relations $S\subseteq X\times Y$ and $T\subseteq Y\times X$ such that $R= ST$.
Otherwise, $R$ is called {\it essential\/}.
The set of all inessential relations on~$X$ span a two-sided ideal~$I$ of~$\CR$. We define $\CE=\CR/I$.
It is clear that $\CE$ is a $k$-algebra having as a $k$-basis the set of all essential relations on~$X$.
The purpose of this paper is to explore the concept of essential relation and to study the structure of~$\CE$.

We shall define a nilpotent ideal $N$ of~$\CE$ and describe completely the quotient $\CP=\CE/N$.
More precisely, $\CP$ is isomorphic to a product of matrix algebras over suitable group algebras,
the product being indexed by the set of all (partial) order relations on~$X$, up to permutation.
Consequently, we know the Jacobson radical $J(\CE)$ and we find all the simple $\CE$-modules.

The idea of passing to the quotient by all elements obtained from something smaller is widely used
in the representation theory of finite groups.
In the theory of $G$-algebras, the notion of Brauer quotient is of this kind (see~\cite{The}).
In the more recent development of the theory of biset functors for finite groups (see~\cite{Bo}),
the same idea plays a key role in~\cite{BST}.
The analogous idea for sets instead of groups yields the notion of essential relation,
which does not seem to have been studied. It is the purpose of this paper to fill this gap.


\vspace{-.1cm}
\section{Essential relations}

\noindent
Given a correspondence $R\subseteq X\times Y$ between a set $X$ and a set~$Y$,
then for every $a\in X$ and $b\in Y$ we write
$$_aR=\{y\in Y \,\mid\, (a,y)\in R \} \qquad \text{and}\qquad R_b=\{x\in X \,\mid\, (x,b)\in R \} \,.$$
We call $_aR$ a {\it column\/} of~$R$ and $R_b$ a {\it row\/} of~$R$.

We first characterize inessential relations. Any subset of $X\times Y$ of the form $U\times V$
will be called a {\it block} (where $U\subseteq X$ and $V\subseteq Y$).

\result{Lemma} \label{blocks} Let $X,Y,Z$ be finite sets.
\begin{enumerate}
\item[(a)] Let $R\subseteq X\times Z$ be a correspondence between $X$ and~$Z$.
Then $R$ factorizes through~$Y$ if and only if $R$ can be decomposed
as a union of blocks indexed by the set~$Y$.
\item[(b)] Let $R$ be a relation on~$X$, where $X$ has cardinality~$n$.
Then $R$ is inessential if and only if $R$ can be decomposed as a union of at most~$n-1$ blocks.
\end{enumerate}
\fresult

\pf
(a) If $R$ factorizes through~$Y$, then $R=ST$, where $S\subseteq X\times Y$ and $T\subseteq Y\times Z$.
Then we can write
$$R=\bigcup_{y\in Y} S_y \times {_yT} \,,$$
as required.

Suppose conversely that $R=\bigcup_{y\in Y}U_y\times V_y$, where each $U_y$ is a subset of~$X$
and each $V_y$ is a subset of~$Z$. Then we define
$$S = \bigcup_{y\in Y} \; U_y\times \{y\} \subseteq X\times Y
\qquad \text{and}\qquad
T=\bigcup_{y\in Y} \; \{y\}\times V_y \subseteq Y\times Z \,.$$
Then $S$ is a correspondence between $X$ and~$Y$, and $U_y=S_y$, the $y$-th row of~$S$.
Similarly, $T$ is a correspondence between $Y$ and~$Z$, and $V_y={_yT}$, the $y$-th column of~$T$.

We now claim that $R=ST$.
If $(x,z)\in R$, then there exists $y\in Y$ such that $(x,z)\in U_y\times V_y$.
It follows that $(x,y)\in S$ and $(y,z)\in T$, hence $(x,z)\in ST$, proving that $R\subseteq ST$.
If now $(x,z)\in ST$, then there exists $y\in Y$ such that $(x,y)\in S$ and $(y,z)\in T$.
It follows that $x\in S_y=U_y$ and $z\in {_yT}=V_y$, hence $(x,z)\in U_y\times V_y \subseteq R$,
proving that $ST\subseteq R$. We have shown that $R=ST$, proving the claim.

(b) This follows immediately from (a).
\endpf

\result{Corollary} \label{equal-rows}
Let $R$ be a relation on~$X$.
If two rows of~$R$ are equal, then $R$ is inessential.
If two columns of~$R$ are equal, then $R$ is inessential.
\fresult

\pf
Suppose that $_aR={_bR}=V$. Then
$$R=(\{a,b\}\times V) \cup
\Big(\bigcup_{{\scriptstyle c\in X}\atop{\scriptstyle c\neq a, c\neq b}}
\{c\}\times {_cR} \, \Big) \,,$$
a union of $n-1$ blocks, where $n=\Card(X)$. The proof for rows is similar.
\endpf

\result{Corollary} \label{empty-row}
Let $R$ be a relation on~$X$.
If a row of~$R$ is empty, then $R$ is inessential.
If a column of~$R$ is empty, then $R$ is inessential.
\fresult

\pf
Assume that $_aR=\emptyset$. Then
$$R=\bigcup_{{\scriptstyle c\in X}\atop{\scriptstyle c\neq a}}
\{c\}\times {_cR} \,,$$
a union of $n-1$ blocks, where $n=\Card(X)$. The proof for rows is similar.
\endpf

\result{Corollary} \label{equivalence}
Let $R$ be a relation on~$X$.
If $R$ is an equivalence relation different from the equality relation
(i.e. $R\neq\Delta$ where $\Delta$ is the diagonal of $X\times X$), then $R$ is inessential.
\fresult

\pf
Suppose that $a$ and $b$ are equivalent and $a\neq b$.
Then the rows $R_a$ and $R_b$ are equal and Corollary~\ref{equal-rows} applies.
\endpf

We need a few basic facts about reflexive relations.
Recall that a relation~$S$ on~$X$ is reflexive if $S$ contains $\Delta=\{(x,x)\mid x\in X\}$.
Moreover, a {\it preorder\/} is a relation which is reflexive and transitive,
while an {\it order\/} is a preorder which is moreover antisymmetric.
(Note that, throughout this paper, the word ``order'' stands for ``partial order''.)
Associated to a preorder~$R$, there is an equivalence relation~$\sim_R$ defined by
$x\sim_R y$ if and only if $(x,y)\in R$ and $(y,x)\in R$.
Then $\sim_R$ is the equality relation if and only if $R$ is an order.

We will often use the containment of relations on~$X$ (as subsets of $X\times X$).
Note that if $R\subseteq S$, then $RT\subseteq ST$ and $TR\subseteq TS$ for any relation~$T$ on~$X$.
If $S$ is a reflexive relation, then $S\subseteq S^2\subseteq S^3\ldots$
and there exists $m\in\N$ such that $S^m=S^{m+1}$, hence $S^m=S^N$ for all $N\geq m$.
The relation $\overline S:=S^m$ is called the {\it transitive closure} of~$S$.
It is reflexive and transitive, that is, a preorder. Note that any preorder is an idempotent relation.

\result{Proposition} \label{preorder}
Let $R$ be a preorder on a finite set~$X$ of cardinality~$n$.
\begin{enumerate}
\item[(a)] If $R$ is not an order, then $R$ is inessential.
\item[(b)] If $R$ is an order and if $Q$ is a reflexive relation contained in~$R$, then $Q$ is essential.
In particular, if $R$ is an order, then $R$ is essential.
\item[(c)] If $R$ is a total order, then $R$ is maximal among essential relations.
\end{enumerate}
\fresult

\pf
(a)  If $R$ is not an order, then the associated equivalence relation $\sim_R$ is not the equality relation.
Let $a$ and $b$ be equivalent under~$\sim_R$ with $a\neq b$.
Then, by transitivity of~$R$, the rows $R_a$ and $R_b$ are equal.
By Corollary~\ref{equal-rows}, $R$ is inessential.

(b) Suppose now that $R$ is an order and that $Q$ is reflexive with $Q\subseteq R$.
We claim that, if $a\neq b$, then $(a,a)$ and $(b,b)$ cannot belong to a block contained in~$Q$.
This is because, if $(a,a),(b,b)\in U\times V \subseteq Q$,
then $(a,b)\in U\times V$ (because $a\in U$ and $b\in V$) and $(b,a)\in U\times V$
(because $b\in U$ and $a\in V$), and therefore $(a,b)\in Q$ and $(b,a)\in Q$,
hence $(a,b)\in R$ and $(b,a)\in R$, contrary to antisymmetry.
It follows that, in any expression of $Q$ as a union of blocks,
the diagonal elements $(a,a)$ all lie in different blocks, so that the number of blocks is at least~$n$.
This shows that $Q$ is essential.

(c) Without loss of generality, we can assume that the total order~$R$ is the usual total order
on the set $X=\{1,2,\ldots,n\}$, i.e. $(x,y)\in R \Leftrightarrow x\leq y$.
Let $S$ be a relation strictly containing~$R$.
Then $S{-}R\neq\emptyset$ and we choose $(j,i)\in S{-}R$ with $i$ maximal,
and then $j$ maximal among all~$x$ with $(x,i)\in S{-}R$.
In other words, $(j,i)\in S$, but $j>i$ because $(j,i)\notin R$, and moreover
$$(x,y)\in S{-}R \;\Rightarrow\; y\leq i \qquad\text{and}\qquad (x,i)\in S{-}R \;\Rightarrow\; x\leq j \,.$$
If $i=j-1$, then the rows $S_i$ and $S_j$ are equal, so $S$ is inessential by Corollary~\ref{equal-rows}.

Assume now that $j-1>i$. Then we claim that 
$$S=\big(S_i \times \{i,j\}\big) \cup \big(S_{j{-}1} \times \{j{-}1,j\}\big) \cup 
\Big(\bigcup_{k\neq i,j{-}1,j} S_k\times \{k\} \Big) \,.$$

To show that the first block is contained in~$S$, let $x\in S_i$.
Then $x\leq j$ if $(x,i)\in S{-}R$, and also $x\leq j$ if $(x,i)\in R$, i.e. $x\leq i$.
Hence $x\leq j$ in both cases, and therefore $(x,j)\in R\subset S$.
This shows that $S_i \times \{j\} \subseteq S$.

To show that the second block is contained in~$S$, let $x\in S_{j{-}1}$.
Then $(x,j{-}1)$ cannot belong to $S{-}R$, by maximality of~$i$ (because $j{-}1>i$).
Thus $(x,j{-}1)\in R$, that is, $x\leq j{-}1$. Then $x<j$, hence $(x,j)\in S$.
This shows that $S_{j{-}1} \times \{j\} \subseteq S$.

Next we show that $S$ is contained in the union of the blocks above.
This is clear for any $(x,y)\in S$ such that $y\neq j$. Now take $(x,j)\in S$.
By maximality of~$i$ and since $j> i$, we have $(x,j)\in R$, that is, $x\leq j$.
If $x=j$, then $(j,j)\in S_i \times \{i,j\}$ because $(j,i)\in S$, that is, $j\in S_i$.
If $x < j$, then $x\leq j{-}1$, hence $(x,j{-}1)\in S$, that is $x\in S_{j{-}1}$,
and therefore $(x,j)\in S_{j{-}1} \times \{j{-}1,j\}$.

This proves the claim about the block decomposition.
It follows that $S$ is a union of $n-1$ blocks, so $S$ is inessential.
\endpf

\result{Example} \label{antidiagonal} 
{\rm Let $n=\Card(X)$.
Let $\Delta$ be the diagonal of $X\times X$ and let $R=(X\times X)-\Delta$.
It is not difficult to see that $R$ is essential if $2\leq n\leq 4$.
However, for $n\geq 5$, we prove that $R$ is inessential.
Without loss of generality, we can choose $X=\{1,\ldots,n\}$.
If $U\subseteq X$, we write $U^c$ for the complement of~$U$ in~$X$.
Then it is easy to prove that $R$ is equal to
$$\Big(\bigcup_{i=1}^{n-3} \{i,i+3\}^c \times  \{i,i+3\} \Big) \cup
\big(\{n{-}2,n{-}1,n\}^c \times \{n{-}2,n{-}1,n\}\big) \cup \big(\{1,2,3\}^c \times \{1,2,3\}\big) .$$
This is a union of $n-1$ blocks, so $R$ is inessential.
}
\fresult


\vspace{-.1cm}
\section{Permutations}

\noindent
As before, $X$ denotes a finite set. 
We let $\Sigma$ be the symmetric group on~$X$, that is, the group of all permutations of~$X$.
For any $\sigma\in\Sigma$, we define
$$\Delta_\sigma= \{\, (\sigma(x),x) \in X\times X \,\mid\, x\in X \,\} \,.$$
This is actually the graph of the map $\sigma^{-1}$, but the choice is made so that
$\Delta_\sigma\Delta_\tau=\Delta_{\sigma\tau}$ for all $\sigma,\tau\in\Sigma$.
With a slight abuse, we shall often call $\Delta_\sigma$ a permutation.
We also write $\Delta=\Delta_{\Id}$.

The group $\Sigma$ has a left action on the set of all relations,
$\sigma$ acting via left multiplication by~$\Delta_\sigma$.
Similarly, $\Sigma$ also acts on the right on the set of relations.
It is useful to note how multiplication by~$\Delta_\sigma$ behaves. Given any relation $R$ on~$X$,
$$(x,y)\in R \quad\Longleftrightarrow \quad (\sigma(x),y)\in \Delta_\sigma R
\quad\Longleftrightarrow \quad (x,\sigma^{-1}(y))\in R\Delta_\sigma \,.$$

\result{Lemma} \label{action}
Let $R$ be a relation on~$X$ and let $\Delta_\sigma$ be a permutation.
\begin{enumerate}
\item[(a)] $R$ is essential if and only if $\Delta_\sigma R$ is essential.
\item[(b)] $\Delta_\sigma$ is essential.
\item[(c)] The left action of $\Sigma$ on the set of all essential relations is free.
\end{enumerate}
\fresult

\pf
(a) If $R$ factorizes through a set of cardinality smaller than~$\Card(X)$, then so does $\Delta_\sigma R$.
The converse follows similarly using multiplication by $\Delta_{\sigma^{-1}}$.

(b) This follows from (a) by taking $R=\Delta$
(which is essential by Proposition~\ref{preorder} because it is an order).

(c) Suppose that $\Delta_\sigma R=R$ for some $\sigma\neq\Id$. Then
$$(x,y)\in R \quad\Longleftrightarrow \quad (\sigma(x),y)\in R \,,$$
hence $_xR={_{\sigma(x)}R}$.
Since $\sigma\neq\Id$, two columns of~$R$ are equal and so $R$ is inessential by Corollary~\ref{equal-rows}.
Thus if $R$ is essential, $\Delta_\sigma R\neq R$ for all $\sigma\neq\Id$.
\endpf

Our next result will be essential in our analysis of essential relations.

\result{Theorem} \label{contain}
Any essential relation contains a permutation.
\fresult

\bigskip
We shall provide two different proofs. The first is direct, while the second uses a theorem of Philip Hall.
In fact, a relation containing a permutation is called a Hall relation in a paper of Schwarz \cite{Sch},
because of Hall's theorem, so Theorem~\ref{contain} asserts that any essential relation is a Hall relation.

\bigskip
\noindent
{\bf First proof~:}
Let $R$ be a relation on~$X$ and let $n=\Card(X)$.
We have to prove that, if $R$ is essential,
then there exists $\sigma\in\Sigma$ such that $R$ contains $\Delta_{\sigma^{-1}}$,
that is, $R\Delta_\sigma$ contains~$\Delta$ (or in other words $R\Delta_\sigma$ is reflexive).
Let $D_\sigma=R\Delta_\sigma \cap\Delta$ and suppose that $\Card(D_\sigma)<n$, for all $\sigma\in\Sigma$.
Then we have to prove that $R$ is inessential.

We choose $\tau\in\Sigma$ such that $\Card(D_\sigma)\leq\Card(D_\tau)$, for all $\sigma\in\Sigma$.
We let $S=R\Delta_\tau$ and we aim to prove that $S$ is inessential (hence $R$ too by Lemma~\ref{action}).
Note that $D_\tau\subseteq S$ by construction. Define
$$A=\{ a\in X \,\mid\, (a,a)\in D_\tau \} \,, \quad\text{ in other words } \quad
D_\tau = \{ (a,a) \,\mid\, a\in A \} \,.$$
In particular $\Card(A)=\Card(D_\tau)<n$.
By maximality of $D_\tau$, we have the following property~:
$$\Card(S\Delta_\sigma \cap\Delta) \leq \Card(A) \,, \; \forall \; \sigma\in\Sigma \,. \qquad(*) $$

Given $x,y\in X$, define a {\it path\/} from $x$ to $y$ to be a sequence $x_0,x_1,\ldots,x_r$ of elements of~$X$
such that $x_0=x$, $x_r=y$, and $(x_i,x_{i+1})\in S$ for all $i=0,\ldots,r-1$.
We write $x\leadsto y$ to indicate that there is a path from $x$ to~$y$,
and also $x\to y$ whenever $(x,y)\in S$ (path of length~1).
Define $A^c$ to be the complement of~$A$ in~$X$ (so $A^c$ is nonempty by assumption). Define also
$$
\begin{array}{rcl}
A_1 &=& \{ a\in A \,\mid\, \text{ there exists }\, z\in A^c \,\text{ and a path } z \leadsto a \} \,, \\
A_2 &=& \{ a\in A \,\mid\, \text{ there exists }\, z\in A^c \,\text{ and a path } a \leadsto z \} \,.
\end{array}
$$

We claim that there is no path from an element of $A_1$ to an element of~$A_2$.
Suppose by contradiction that there is a path $a_1 \leadsto a_2$ with $a_i\in A_i$.
Then there exists $z_i\in A^c$ and paths
$z_1 \leadsto a_1  \leadsto a_2  \leadsto z_2$, in particular $z_1 \leadsto a \leadsto z_2$ with $a\in A$.
In the path $z_1 \leadsto a$, let $w_1$ be the element of~$A^c$ closest to~$a$,
so that the path $w_1 \leadsto a$ does not contain any element of~$A^c$ except~$w_1$.
Similarly, let $w_2$ be the element of~$A^c$ closest to~$a$ in the path $a\leadsto z_2$ ,
so that the path $a \leadsto w_2$ does not contain any element of~$A^c$ except~$w_2$.
We obtain a path $w_1 \leadsto a \leadsto w_2$ having all its elements in~$A$
except the two extremities $w_1$ and $w_2$.
By suppressing cycles within~$A$, we can assume that all elements of~$A$ in this path are distinct.
We end up with a path
$$w_1\to x_1\to\cdots\to x_r\to w_2$$
where $x_1,\ldots,x_r\in A$ are all distinct.

Let $\sigma\in\Sigma$ be the cycle defined by $\sigma(w_1)=x_1$, $\sigma(x_i)=x_{i+1}$ for $1\leq i\leq r-1$,
$\sigma(x_r)=w_2$, and finally $\sigma(w_2)=w_1$ in case $w_2\neq w_1$.
In case $w_2=w_1$, then $\sigma(w_2)$ is already defined to be $\sigma(w_2)=\sigma(w_1)=x_1$.
We emphasize that $\sigma(y)=y$ for all the other elements $y\in X$. Then we obtain~:
$$
\begin{array}{rcll}
(w_1,x_1)\in S \qquad &\text{hence}& \qquad (w_1,w_1)\in S\Delta_\sigma \,, &\\
(x_i,x_{i+1})\in S \qquad &\text{hence}& \qquad (x_i,x_i)\in S\Delta_\sigma \,, &\\
(x_r,w_2)\in S \qquad &\text{hence}& \qquad (x_r,x_r)\in S\Delta_\sigma \,, &\\
(y,y)\in S \qquad &\text{hence}& \qquad (y,y)\in S\Delta_\sigma \,, &
\forall \, y\in A-\{x_1,\ldots,x_r\} \,.\\
\end{array}
$$
Thus we obtain $(a,a)\in S\Delta_\sigma$, $\forall \, a\in A$, but also $(w_1,w_1)\in S\Delta_\sigma$.
Therefore $\Card(S\Delta_\sigma \cap \Delta)>\Card(A)$, contrary to Property~$(*)$.
This proves the claim that there is no path from $A_1$ to~$A_2$.

In particular, $A_1\cap A_2=\emptyset$ because if $a\in A_1\cap A_2$
we would have a path of length zero from $A_1$ to~$A_2$ (since $(a,a)\in S$).
Let $A_3$ be the complement of $A_1\cup A_2$ in~$A$.
Thus $X$ is the disjoint union of the 4 subsets $A^c$, $A_1$, $A_2$, and $A_3$.

We now claim the following~:
\begin{itemize}
\item[(a)] There is no relation between $A^c$ and $A^c$, that is, $S\cap(A^c\times A^c)=\emptyset$.
\item[(b)] There is no relation between $A^c$ and $A_2\cup A_3$, that is,
$S\cap(A^c\times A_2)=\emptyset$ and $S\cap(A^c\times A_3)=\emptyset$.
\item[(c)] There is no relation between $A_1$ and $A^c$, that is, $S\cap(A_1\times A^c)=\emptyset$.
\item[(d)] There is no relation between $A_1$ and $A_2$, that is, $S\cap(A_1\times A_2)=\emptyset$.
\item[(e)] There is no relation between $A_1$ and $A_3$, that is, $S\cap(A_1\times A_3)=\emptyset$.
\end{itemize}

To prove (a), suppose that $(w,z)\in S$ where $w,z\in A^c$.
Choose $\rho\in\Sigma$ such that $\rho(a)=a$ for all $a\in A$ and $\rho(w)=z$.
Then $(a,a)\in S\Delta_\rho$ for all $a\in A$ but also $(w,w)\in S\Delta_\rho$,
contrary to Property~$(*)$.

To prove (b), we note that the definition of $A_1$ implies that, if $(z,a)\in S$ with $z\in A^c$ and $a\in A$,
then $a\in A_1$. Thus $a\notin A_2\cup A_3$.

To prove (c), suppose that $(a,z)\in S$ with $a\in A$ and $z\in A^c$.
Then $a\in A_2$ by the definition of~$A_2$ and in particular $a\notin A_1$.

Property (d) follows immediately from the previous claim that there is no path from $A_1$ to $A_2$.

To prove (e), suppose that $(a_1,a_3)\in S$ with $a_1\in A_1$ and $a_3\in A_3$.
Then by the definition of $A_1$, there is a path $z\leadsto a_1\to a_3$ where $z\in A^c$,
but this means that $a_3\in A_1$, a contradiction.

It follows that the relation $S$ has the property that there is no relation
between $A^c\cup A_1$ and $A^c\cup A_2\cup A_3$.
Therefore $S$ is the union of the columns indexed by $(A^c\cup A_1)^c=A_2\cup A_3$
and the lines indexed by $(A^c\cup A_2\cup A_3)^c=A_1$, that is,
$$S= \Big(\bigcup_{b\in A_2\cup A_3} {_bS}\Big) \cup \Big(\bigcup_{a\in A_1} S_a\Big) \,.$$
Since $\Card(A_2\cup A_3)+\Card(A_1) =\Card(A)$, we obtain a union of $\Card(A)$ blocks.
But $\Card(A)<n$ by assumption, so $S$ is inessential, as was to be shown.
\endpf

\noindent
{\bf Second proof~:}
Let $R$ be an essential relation on~$X$. For any subset $A$ of~$X$, define
$$R_A=\{ x\in X \mid \exists\, a\in A \text{ such that } (x,a)\in R\} = \bigcup_{a\in A} R_a \,.$$
Then $R$ decomposes as a union of blocks
$$R=\Big( \bigcup_{y\notin A} (R_y\times\{y\}) \Big) \bigcup
\Big( \bigcup_{x\in R_A} (\{x\}\times {_xR}) \Big) \,.$$
Since $R$ is essential, $\Card(X-A)+\Card(R_A)$ cannot be strictly smaller than $\Card(X)$.
Therefore $\Card(R_A)\geq\Card(A)$, for all subsets $A$ of~$X$, that is
$$\Card\big(\bigcup_{a\in A} R_a\big)\geq \Card(A) \,.$$
This is precisely the assumption in a theorem of Philip Hall
(see Theorem 5.1.1 in~\cite{HaM}, or \cite{HaP} for the original version which is slightly different).
The conclusion is that there exist elements $x_y\in R_y$, where $y$ runs over~$X$, which are all distinct.
In other words $\sigma:y\mapsto x_y$ is a permutation and
$$(\sigma(y),y)=(x_y,y) \in R \quad\text{ for all } y\in X \,.$$
This means that $R$ contains $\Delta_\sigma$, as required.
\endpf

\result{Corollary} \label{power}
Let $R$ be an essential relation on~$X$. Then there exists $m\in \N$ such that $R^m$ is a preorder.
\fresult

\pf
By Theorem~\ref{contain}, $R$ contains $\Delta_\sigma$ for some $\sigma\in\Sigma$.
If $\sigma$ has order~$k$ in the group~$\Sigma$, then $R^k$ contains $\Delta_{\sigma^k}=\Delta$,
so $R^k$ is reflexive. Then the transitive closure of~$R^k$ is some power $R^{kt}$.
This is reflexive and transitive, that is, a preorder.
\endpf

We know that any order is an essential relation (Proposition~\ref{preorder}),
hence contains a permutation (Theorem~\ref{contain}). But in fact, we have a more precise result.

\result{Lemma} \label{unique}
If $R$ is an order on $X$, then $R$ contains a unique permutation, namely~$\Delta$.
\fresult

\pf
Suppose that $R$ is reflexive and transitive and contains a nontrivial permutation $\Delta_\sigma$.
Then $\sigma$ contains a nontrivial $k$-cycle, say on $x_1,\ldots,x_k$, for some $k\geq2$.
It follows that $(x_{i+1},x_i)\in R$ for $1\leq i\leq r-1$, hence $(x_k,x_1)\in R$ by transitivity of~$R$.
Now we also have $(x_1,x_k)\in R$ because $\sigma(x_k)=x_1$.
Thus the relation $R$ is not antisymmetric, hence cannot be an order.
\endpf

In the same vein, we have the following more general result.

\result{Lemma} \label{S'S}
Let $R$ be an order and let $S,S'$ be two relations on $X$. The following two conditions are equivalent~:
\begin{itemize}
\item[(a)] $\Delta\subseteq S'S \subseteq R$.
\item[(b)] There exists a permutation $\Delta_\sigma$ such that~:

$\Delta\subseteq \Delta_{\sigma^{-1}}S\subseteq R\;$ and $\;\Delta\subseteq S'\Delta_\sigma\subseteq R$.
\end{itemize}
Moreover, in condition~(b), the permutation $\sigma$ is unique.
\fresult

\pf
If (b) holds, then
$$\Delta=\Delta^2\subseteq (S'\Delta_\sigma)(\Delta_{\sigma^{-1}}S)=S'S \subseteq R^2=R \,,$$
so (a) holds.

If (a) holds, then $S'S$ is essential, by Proposition~\ref{preorder}.
It follows that $S$ is essential, and therefore $S$ contains a permutation $\Delta_\sigma$,
by Theorem~\ref{contain}. Then we obtain
$$\Delta\subseteq \Delta_{\sigma^{-1}}S \qquad \text{ and } \qquad
S'\Delta_\sigma\subseteq S'S\subseteq R \,.$$
Similarly, $S'$ is essential, hence contains a permutation $\Delta_\tau$, and we obtain
$$\Delta\subseteq S'\Delta_{\tau^{-1}} \qquad \text{ and } \qquad
\Delta_\tau S\subseteq S'S\subseteq R \,.$$
Now $R$ contains $\Delta_\tau\Delta_\sigma=\Delta_{\tau\sigma}$ and Lemma~\ref{unique} implies that
$\tau\sigma=\Id$, that is, $\tau=\sigma^{-1}$. Then (b) follows.

If moreover $S$ contains a permutation $\Delta_\rho$ (that is, $\Delta\subseteq \Delta_{\rho^{-1}}S$), then
$$\Delta_{\sigma^{-1}\rho}=\Delta_{\sigma^{-1}}\Delta_\rho\subseteq S'S\subseteq R \,,$$
and so $\sigma^{-1}\rho=\Id$ by Lemma~\ref{unique}, proving the uniqueness of~$\sigma$.
\endpf


\vspace{-.1cm}
\section{The essential algebra}

\noindent
Let $X$ be a finite set and let $k$ be a commutative ring.
We shall be mainly interested in the cases where $k$ is either the ring~$\Z$ of integers or a field,
but it is convenient to work with an arbitrary commutative ring.

Let $\CR$ be the $k$-algebra of the monoid of all relations on~$X$.
This monoid is a $k$-basis of~$\CR$ and the product in the monoid defines the algebra structure.
The set of all inessential relations on~$X$ spans a two-sided ideal~$I$ of~$\CR$.
We define $\CE=\CR/I$ and call it the {\it essential algebra\/}.
It is clear that $\CE$ is a $k$-algebra having as a $k$-basis the set of all essential relations on~$X$.
Moreover, if $R$ and $S$ are essential relations but $RS$ is inessential, then $RS=0$ in~$\CE$.

Both $\CR$ and $\CE$ have an anti-automorphism, defined on the basis elements by $R\mapsto R^{op}$,
where $(x,y)\in R^{op}$ if and only if $(y,x)\in R$. It is easy to see that $(RS)^{op}=S^{op}R^{op}$.

As before, we let $\Sigma$ be the symmetric group of all permutations of~$X$.
We first describe an obvious quotient of~$\CE$.

\result{Lemma} \label{group-algebra}
Let $H$ be the $k$-submodule of the essential algebra~$\CE$ spanned by
the set of all essential relations which strictly contain a permutation.
Then $H$ is a two-sided ideal of~$\CE$ and $\CE/H\cong k\Sigma$,
the group algebra of the symmetric group~$\Sigma$.
\fresult

\pf
Let us write $\subset$ for the strict containment relation.
Let $R$ be an essential relation such that $\Delta_\sigma\subset R$ and let $S$ be any essential relation.
Then $S$ contains a permutation $\Delta_\tau$, by Theorem~\ref{contain}. We obtain
$$\Delta_{\sigma\tau}=\Delta_\sigma\Delta_\tau \subset R\Delta_\tau\subseteq RS \,,$$
showing that $RS\in H$. Similarly $SR\in H$ and therefore $H$ is a two-sided ideal of~$\CE$.

The quotient $\CE/H$ has a $k$-basis consisting of all the permutations $\Delta_\sigma$, for $\sigma\in\Sigma$.
Moreover, they multiply in the same way as permutations,
so $\CE/H$ is isomorphic to the group algebra of the symmetric group~$\Sigma$.
\endpf

If $k$ is a field, it follows, not surprisingly,
that every irreducible representation of the symmetric group~$\Sigma$ gives rise to a simple $\CE$-module.
In short, the representation theory of the symmetric group~$\Sigma$ is part of  the representation theory of~$\CE$.

We now want to describe another $\CE$-module, which is simple when $k$ is a field.
We fix a total order $T$ on~$X$ (e.g. the usual total order on $X=\{1,\ldots,n\}$).
Then any other total order on~$X$ is obtained by permuting the elements of~$X$.
Since permuting via $\sigma$ corresponds to conjugation by~$\Delta_\sigma$, we see that
$\{T_\sigma:=\Delta_\sigma T \Delta_{\sigma^{-1}} \mid \sigma\in \Sigma \}$
is the set of all total orders on~$X$.
All of them are maximal essential relations on~$X$, by Proposition~\ref{preorder}.

\result{Lemma} \label{total} Let $T$ be a total order on~$X$.
\begin{itemize}
\item[(a)] If $\rho\in\Sigma$, then $T\Delta_\rho T=0$ in~$\CE$ if $\rho\neq\Id$ and otherwise $T\Delta T=T^2=T$.
\item[(b)] The set $\{T_\sigma:=\Delta_\sigma T \Delta_{\sigma^{-1}} \mid \sigma\in \Sigma \}$
is a set of pairwise orthogonal idempotents of~$\CE$.
\end{itemize}
\fresult

\pf
(a) $TT_\rho$ contains both $T$ and $T_\rho$ (because both $T$ and $T_\rho$ contain~$\Delta$).
Since $T\neq T_\rho$ if $\rho\neq\Id$, the product $TT_\rho$ contains strictly $T$
and is therefore inessential by Proposition~\ref{preorder}.
Thus $T\Delta_\rho T$ is also inessential, that is, $T\Delta_\rho T=0$.
On the other hand $T^2=T$ because any preorder is idempotent.

(b) It follows from (a) that
$$T_\sigma T_\tau =\Delta_\sigma T \Delta_{\sigma^{-1}\tau} T \Delta_{\tau^{-1}}=
\left\{\begin{array}{ll}
0    \quad   &\text{if }\; \sigma\neq\tau \,, \\
\Delta_\sigma T\Delta_{\sigma^{-1}}=T_\sigma \quad &\text{if }\; \sigma=\tau \,,
\end{array}\right.
$$
as was to be shown.
\endpf

\result{Proposition} \label{n!}
Fix a total order $T$ on~$X$.
Let $L$ be the $k$-submodule of the essential algebra~$\CE$ spanned
by the set $\{\Delta_\sigma T\mid  \sigma\in \Sigma \}$.
Then $L$ is a left ideal of~$\CE$ and is free of rank~$n!$~as a $k$-module, where $n=\Card(X)$.
If $k$ is a field, then $L$ is a simple $\CE$-module of dimension~$n!$~.
\fresult

\pf
Write $S_\sigma=\Delta_\sigma T$, for all $\sigma\in\Sigma$.
Let $R$ be an essential relation on~$X$.
Then $\Delta_\tau\subseteq R$ for some $\tau\in \Sigma$ by Theorem~\ref{contain}.
Therefore $\Delta\subseteq\Delta_{\tau^{-1}} R$ and this implies that
$$S_\sigma=\Delta S_\sigma\subseteq \Delta_{\tau^{-1}} RS_\sigma 
\qquad\text{and}\qquad
T= \Delta_{\sigma^{-1}} S_\sigma\subseteq \Delta_{\sigma^{-1}}\Delta_{\tau^{-1}} RS_\sigma\,.$$
If this containment is strict, then $\Delta_{\sigma^{-1}}\Delta_{\tau^{-1}} RS_\sigma$ is inessential
(by Proposition~\ref{preorder}) and so $RS_\sigma$ is inessential too (by Lemma~\ref{action}).
Otherwise $S_\sigma=\Delta_{\tau^{-1}} RS_\sigma$, hence $RS_\sigma=\Delta_\tau S_\sigma
=\Delta_{\tau\sigma} T=S_{\tau\sigma}$.
Therefore, in the algebra~$\CE$, either $RS_\sigma=0$ or $RS_\sigma=S_{\tau\sigma}$.
This proves that $L$ is a left ideal of~$\CE$.

Clearly $L$ has rank~$n!$ with basis $\{S_\sigma \mid  \sigma\in \Sigma \}$.
The action of $\CE$ on~$L$ induces a $k$-algebra map
$$\phi: \CE \longrightarrow M_{n!}(k)$$
and $L$ can be viewed as an $M_{n!}(k)$-module
(consisting of column vectors with entries in~$k$).
By Lemma~\ref{total}, the action of $\Delta_\tau T\Delta_{\rho^{-1}}$ on basis elements is given by
$$(\Delta_\tau T\Delta_{\rho^{-1}})\cdot S_\sigma
=(\Delta_\tau T\Delta_{\rho^{-1}})\cdot \Delta_\sigma T
=\left\{\begin{array}{ll}
0    \quad   &\text{if }\; \rho\neq\sigma \,, \\
S_\tau  \quad &\text{if }\; \rho=\sigma \,.
\end{array}\right.
$$
This means that $\phi(\Delta_\tau T\Delta_{\rho^{-1}})$ is the elementary matrix
with a single nonzero entry~1 in position $(\tau,\rho)$. Therefore the map $\phi$ is surjective.
This implies that, if $k$ is a field, the module $L$ is simple as an $\CE$-module, because
the space of column vectors is a simple $M_{n!}(k)$-module.
\endpf


\vspace{-.1cm}
\section{A nilpotent ideal}

\noindent
The purpose of this section is to construct a suitable nilpotent ideal~$N$ of the essential algebra~$\CE$.
We shall later pass to the quotient by~$N$ and describe the quotient $\CE/N$.
In order to find nilpotent ideals, the following well-known result is often useful.

\result{Lemma} \label{nilpotent}
Let $k$ be a commutative ring and let $\CB$ be a $k$-algebra which is finitely generated as a $k$-module.
Let $I$ be a two-sided ideal of~$\CB$ which is $k$-linearly spanned by a set of nilpotent elements of~$\CB$.
\begin{itemize}
\item[(a)] If $k$ is a field, then $I$ is a nilpotent ideal of~$\CB$.
\item[(b)] If $k=\Z$ and if $\CB$ is a finitely generated free $\Z$-module, then $I$ is a nilpotent ideal of~$\CB$.
\item[(c)] Suppose that $\CB$ is defined over~$\Z$, that is, $\CB\cong k\otimes_{\Z} \CB_{\Z}$
for some $\Z$-algebra~$\CB_{\Z}$ which is finitely generated free as a $\Z$-module.
Suppose that $I$ is defined over~$\Z$, that is, $I\cong k\otimes_{\Z} I_{\Z}$,
where $I_{\Z}$ is a two-sided ideal of~$\CB_{\Z}$ which is $\Z$-linearly spanned
by a set of nilpotent elements of~$\CB_{\Z}$.
Then $I$ is a nilpotent ideal of~$\CB$.
\end{itemize}
\fresult

\pf
(a) The assumption still holds after extending scalars to an algebraic closure of~$k$.
Therefore we can assume that $k$ is algebraically closed. Let $J(\CB)$ be the Jacobson radical of~$\CB$.
Then $J(\CB)=\bigcap_{i=1}^rM_i$, where $M_i$ is a maximal two-sided ideal of~$\CB$.
Moreover, by Wedderburn's theorem, $\CB/M_i$ is isomorphic to a matrix algebra $M_{n_i}(k)$,
because $k$ is algebraically closed. We will show that $I\subseteq M_i$, for all $i=1,\ldots,r$.
It then follows that $I\subseteq J(\CB)$, so $I$ is nilpotent
(because it is well-known that the Jacobson radical of a finite-dimensional $k$-algebra is nilpotent).

Let $\overline I$ be the image of~$I$ in $\CB/M_i$.
Then $\overline I$ is spanned by nilpotent elements of~$M_{n_i}(k)$.
But any nilpotent matrix has trace zero
(because its characteristic polynomial is $X^{n_i}$ and the coefficient of $X^{n_i-1}$ is the trace, up to sign).
It follows that $\overline I$ is contained in $\Ker(\tr)$, which is a proper subspace of~$M_{n_i}(k)$.
Now $\overline I$ is a two-sided proper ideal of the simple algebra $M_{n_i}(k)$,
hence $\overline I=\{0\}$, proving that $I\subseteq M_i$.

(b) Let $F$ be a basis of $\CB$ as a $\Z$-module.
Extending scalars to~$\Q$, we see that $F$ is a $\Q$-basis of the $\Q$-algebra $\Q\otimes_{\Z}\CB$
and $\CB$ embeds in~$\Q\otimes_{\Z}\CB$.
By part~(a), the ideal $\Q\otimes_{\Z}I$ is nilpotent in $\Q\otimes_{\Z}\CB$.
Since $I$ embeds in~$\Q\otimes_{\Z}I$, it follows that $I$ is nilpotent.

(c) By part~(b), $I_{\Z}$ is a nilpotent ideal of~$\CB_{\Z}$.
Extending scalars to~$k$, we see that $I$ is a nilpotent ideal of~$\CB$.
\endpf

Recall that $\Sigma$ denotes the symmetric group on~$X$
and that, if $R$ is a reflexive relation, then $\overline R$ denotes the transitive closure of~$R$.

\result{Lemma} \label{conj}
If $S=\Delta_{\tau^{-1}} R\Delta_\tau$ where $\tau\in\Sigma$ and $R$ is a reflexive relation, then
$\overline S=\Delta_{\tau^{-1}} \overline R\Delta_\tau$.
\fresult

\pf
We have $\overline R=R^m$ for some $m$ and we obtain
$$S^m=(\Delta_{\tau^{-1}} R\Delta_\tau)^m=\Delta_{\tau^{-1}} R^m\Delta_\tau
=\Delta_{\tau^{-1}} \overline R\Delta_\tau \,.$$
Therefore $S^m$ is a preorder, because it is conjugate to a preorder, and so
$\overline S=S^m=\Delta_{\tau^{-1}} \overline R\Delta_\tau$.
\endpf

\result{Theorem} \label{ideal}
Let $\CF$ be the set of all reflexive essential relations on~$X$.
Let $N$ be the $k$-submodule of the essential algebra~$\CE$ generated by all elements of the form
$(S-\overline S)\Delta_\sigma$ with $S\in\CF$ and $\sigma\in\Sigma$
(where $\overline S$ denotes the transitive closure of~$S$).
\begin{itemize}
\item[(a)] $N$ is a nilpotent two-sided ideal of $\CE$.
In particular, $N$ is contained in the Jacobson radical~$J(\CE)$.
\item[(b)] The quotient algebra $\CP=\CE/N$ has a $k$-basis consisting of all elements of the form
$S\Delta_\sigma$, where $S$ runs over the set of all orders on~$X$
and $\sigma$ runs over the symmetric group~$\Sigma$.
\end{itemize}
\fresult

\pf
(a) Let $\CE_1$ be the subalgebra of~$\CE$ which is $k$-linearly generated
by the set~$\CF$ of all reflexive essential relations.
It is clearly a subalgebra since the product of two reflexive relations is reflexive.
Let $N_1$ be the $k$-submodule of~$\CE_1$ generated by all elements of the form
$S-\overline S$ with $S\in\CF$. We claim that $N_1$ is a two-sided ideal of~$\CE_1$.

If $T\in\CF$, then $\overline{T\overline{S}}=\overline{TS}$ (because $\overline{TS}$ contains both $T$
and $\overline{S}$, hence $\overline{T\overline{S}}$, and $\overline{T\overline{S}}$ contains both $T$
and $S$, hence $\overline{TS}$). Therefore
$$T(S-\overline S)=(TS-\overline{TS})-(T\overline S-\overline{TS})
=(TS-\overline{TS})-(T\overline S-\overline{T\overline{S}})\,.$$
Note that if $TS$ is inessential (hence zero in~$\CE_1$), then its transitive closure
$\overline{TS}$ cannot be an order by Proposition~\ref{preorder} and is therefore also zero in~$\CE_1$
(again by Proposition~\ref{preorder}).
Thus, in the expression above, we obtain either generators of~$N_1$ or zero.
The same argument works for right multiplication by~$T$
(or use the anti-automorphism of~$\CE$) and this proves the claim.

The ideal $N_1$ is invariant under conjugation by~$\Sigma$ because, for every $\sigma\in\Sigma$,
$$\Delta_{\sigma^{-1}} (S-\overline S)\Delta_\sigma
=\Delta_{\sigma^{-1}} S\Delta_\sigma - \Delta_{\sigma^{-1}} \overline S\Delta_\sigma
=\Delta_{\sigma^{-1}} S\Delta_\sigma - \overline{\Delta_{\sigma^{-1}} S\Delta_\sigma}$$
by Lemma~\ref{conj}.
Therefore the generators of~$N$ can also be written $\Delta_\sigma(S'-\overline{S'})$
with $S'\in\CF$ and $\sigma\in\Sigma$ (namely $S'=\Delta_{\sigma^{-1}} S \Delta_\sigma$).
It follows that $N=N_1\Delta_\Sigma=\Delta_\Sigma N_1$, where we write for simplicity
$\Delta_\Sigma=\{\Delta_\sigma \mid \sigma\in\Sigma \}$.

If $R$ is an essential relation on~$X$, then $R$ contains a permutation $\Delta_\sigma$
(for some $\sigma\in\Sigma$) by Theorem~\ref{contain}, so $R=Q\Delta_\sigma$ with $Q\in\CF$,
and also $R=\Delta_\sigma Q'$ where $Q'=\Delta_{\sigma^{-1}} Q \Delta_\sigma$.
Since $N_1$ is an ideal of~$\CE_1$, it follows that $N$ is invariant by right and left multiplication by~$R$.
Thus $N$ is a two-sided ideal of~$\CE$.

The generators of $N_1$ are nilpotent, because if $\overline S=S^m$, then
$$\begin{array}{rl}
(S-\overline{S})^m=(S-S^m)^m=&
\displaystyle\sum_{j=0}^m {m\choose j} (-1)^j S^{m-j}S^{mj} \\
=&\displaystyle\Big(\sum_{j=0}^m {m\choose j} (-1)^j \Big) S^m = (1-1)^m S^m=0 \,.
\end{array}
$$
Thus $N_1$ is a nilpotent ideal of~$\CE_1$, by Lemma~\ref{nilpotent}
(because clearly $\CE_1$ and $N_1$ are defined over~$\Z$).
Since $N_1$ is invariant under conjugation by~$\Sigma$, we obtain
$N^n=(N_1\Delta_\Sigma)^n=N_1^n\Delta_\Sigma$ for every $n\in\N$.
Since $N_1^m=0$ for some~$m$, the ideal $N$ is nilpotent.

(b) In the quotient algebra $\CE/N$, any reflexive relation $Q$ is identified with its transitive closure~$\overline Q$.
Moreover, by Theorem~\ref{contain}, any essential relation $R$ on~$X$ can be written
$R=Q\Delta_\sigma$, with $Q$ reflexive,
and $Q\Delta_\sigma$ is identified with $\overline Q\Delta_\sigma$ in the quotient algebra~$\CE/N$.
Note that $\overline Q$ is a preorder and that $\overline Q$ is zero in~$\CE$ if it is not an order,
by Proposition~\ref{preorder}.

On each basis element~$R$ of~$\CE$, the effect of passing to the quotient by~$N$
consists of just two possibilities.
\begin{itemize}
\item[$\bullet$] If $R$ can be written $R=Q\Delta_\sigma$, with $Q$ reflexive and $\overline Q$ is not an order,
then $R$ is identified with $\overline Q\Delta_\sigma$, so $R$ is zero in~$\CE/N$ because $\overline Q$ is zero.
\item[$\bullet$] If $R$ can be written $R=Q\Delta_\sigma$, with $Q$ reflexive and $\overline Q$ is an order,
then $R$ is identified with an element of the form $S\Delta_\sigma$
where $S$ is an order (namely $S=\overline Q$).
\end{itemize}
In the second case, $\Delta_\sigma$ is the unique permutation contained in~$R$
(or in other words the expression $R=Q\Delta_\sigma$ is the unique decomposition of $R$ as a product
of a reflexive relation and a permutation). This is because if 
$\Delta_{\sigma'}\subseteq R$, we obtain
$$\Delta_{\sigma'\sigma^{-1}}=\Delta_{\sigma'}\Delta_{\sigma^{-1}}\subseteq R\Delta_{\sigma^{-1}}=Q
\subseteq \overline Q \,,$$
so that $\sigma'\sigma^{-1}=\Id$ since $\overline Q$ is an order (Lemma~\ref{unique}).
Thus $\sigma'=\sigma$. This uniqueness property shows, on the one hand,
that both possibilities cannot occur simultaneously and, on the other hand,
that in the second case the order $S=\overline Q$ is uniquely determined by~$R$.

It follows that the nonzero images in~$\CE/N$ of the basis elements of~$\CE$ form a $k$-basis of~$\CE/N$
consisting of (the images of) the elements $S\Delta_\sigma$, where $S$ is an order and $\sigma\in\Sigma$.\endpf

The quotient algebra $\CP=\CE/N$ will be called the {\it algebra of permuted orders\/} on~$X$,
because every basis element $S\Delta_\sigma$ is obtained from the order~$S$
by applying a permutation~$\sigma$ to the rows of~$S$.
Moreover, $\Delta_\sigma$ is the unique permutation contained in~$S\Delta_\sigma$,
because $\Delta$ is the unique permutation contained in~$S$ by Lemma~\ref{unique}.
This defines a $\Sigma$-grading on~$\CP$~:
$$\CP = \bigoplus_{\sigma\in \Sigma} \CP_\sigma \,,$$
where $\CP_\sigma$ is spanned by the set of all permuted orders containing~$\Delta_\sigma$.
Clearly $\CP_\sigma \cdot \CP_\tau = \CP_{\sigma\tau}$, so we have indeed a $\Sigma$-grading.
We also write $\CP_1:=\CP_{\Id}$ and call it the {\it algebra of orders\/} on~$X$.
Moreover, $\CP_\sigma=\Delta_\sigma\CP_1=\CP_1\Delta_\sigma$,
so that the product in~$\CP$ is completely determined by the product in the subalgebra~$\CP_1$
and the product in the symmetric group~$\Sigma$.
Hence we first need to understand the subalgebra~$\CP_1$.


\vspace{-.1cm}
\section{The algebra of orders}

\noindent
Let $\CP_1$ be the algebra of orders on~$X$ defined above.
It has a $k$-basis $\CO$ consisting of all orders on~$X$.
The product of basis elements $R,S\in\CO$ will be written $R\cdot S$ and is described as follows.

\result{Lemma} \label{order-product} Let $\cdot$ be the product in the $k$-algebra $\CP_1$.
\begin{itemize}
\item[(a)] Let $R,S\in\CO$. Then the product $R\cdot S$ is equal to the transitive closure of $R\cup S$
if this closure is an order, and zero otherwise.
\item[(b)] The product $\cdot$ is commutative.
\end{itemize}
\fresult

\pf
(a) By definition of the ideal $N$, the product $RS$ in the algebra $\CP=\CE/N$ is identified with
the transitive closure $\overline{RS}$, which is also the transitive closure of $R\cup S$,
because the inclusions
$$R\cup S\subseteq RS \subseteq (R\cup S)^2\subseteq \overline{R\cup S}\subseteq\overline{RS}$$
force the equality $\overline{R\cup S}=\overline{RS}$.
Now $\overline{RS}$ is a preorder. If this is an order, then $R\cdot S=\overline{RS}$.
If this preorder is not an order, then it is zero in~$\CE$ (by Proposition~\ref{preorder}),
hence also zero in~$\CP_1$.

(b) This follows from (a) and the fact that $R\cup S=S\cup R$.
\endpf

\result{Theorem} \label{semi-simple}
\begin{itemize}
\item[(a)] There exists a $k$-basis $\{f_R \,|\, R\in\CO\}$ of $\CP_1$,
consisting of mutually orthogonal idempotents whose sum is~$1$,
and such that, for every $R\in\CO$, the ideal generated by $f_R$ is free of rank one as a $k$-module.
\item[(b)] $\CP_1$ is isomorphic to a product of copies of~$k$, indexed by~$\CO$~:
$$\CP_1 \cong \prod_{R\in\CO} k{\cdot} f_R \,.$$
\end{itemize}

\fresult

\pf
We know that $\CP_1$ is commutative, with a basis~$\CO$ consisting of all orders on~$X$.
Any such basis element is idempotent.
Moreover $\CO$ is a partially ordered set with respect to the containment relation
and we make it a lattice by adding an element $\infty$ and defining
$R\vee S=\infty$ whenever the transitive closure of $R\cup S$ is not an order,
while $R\vee S$ is the transitive closure of $R\cup S$ otherwise.
The greatest lower bound of $R$ and $S$ is just the intersection $R\cap S$.

Now define $g_R=R$ if $R\in\CO$ and $g_\infty=0$.
By Lemma~\ref{order-product}, these elements satisfy the condition $g_R\cdot g_S=g_{R\vee S}$.
Therefore Theorem~\ref{inversion} of the appendix applies. We let
$$f_R=\sum_{\scriptstyle S\in\CO \atop\scriptstyle R\subseteq S} \mu(R,S)S \,,$$
where $\mu$ denotes the M\"obius function of the poset~$\CO$, so by M\"obius inversion, we have
$$R=\sum_{\scriptstyle S\in\CO \atop\scriptstyle R\subseteq S} f_S \,.$$
The transition matrix from $\{R\in\CO\}$ to $\{f_R\mid R\in\CO\}$ is upper-triangular,
with 1 along the main diagonal, hence invertible over~$\Z$.
It follows that $\{f_R\mid R\in\CO\}$ is a $k$-basis of~$\CP_1$.
By Theorem~\ref{inversion} of the appendix, $\{f_R \,|\, R\in\CO\}$
is a set of mutually orthogonal idempotents in~$\CP_1$ whose sum is~$1$.
Moreover, by the same theorem,
$$f_R\cdot T=
\left\{\begin{array}{ll}
f_R\quad &\text{ if } T\subseteq R \,, \\
\,0 \quad &\text{ if } T\not\subseteq R \,. \\
\end{array}\right.$$
Since $T$ runs over a $k$-basis of $\CP_1$,
this proves that the ideal $\CP_1f_R$ generated by $f_R$ is equal to
the rank one submodule $k{\cdot}f_R$ spanned by~$f_R$.
Thus we obtain
$$\CP_1 \cong \prod_{R\in\CO} \CP_1f_R = \prod_{R\in\CO} k{\cdot} f_R \,,$$
as was to be shown.
\endpf

Note that if $k$ is a field, then each idempotent $f_R$ is primitive.


\vspace{-.1cm}
\section{The algebra of permuted orders}

\noindent
We know from the end Section~5 that the algebra $\CP$ of permuted orders is $\Sigma$-graded
$$\CP=\bigoplus_{\sigma\in\Sigma}\CP_\sigma \,.$$
If $R,S\in\CO$ and $\sigma,\tau\in\Sigma$, then the product in~$\CP$ satisfies
$$(R\Delta_\sigma)(S\Delta_\tau)
= \big(R\cdot(\Delta_\sigma S\Delta_{\sigma^{-1}})\big)\Delta_\sigma\Delta_\tau
= \big(R\cdot(\Delta_\sigma S\Delta_{\sigma^{-1}})\big)\Delta_{\sigma\tau} \,,$$
where $\cdot$ denotes the product in~$\CP_1$ described in Lemma~\ref{order-product} .
Note that this definition makes sense because $\Delta_\sigma S\Delta_{\sigma^{-1}}$ is an order,
since $S$ is.
Note also that we can write the basis elements as $\Delta_\sigma S$ with $S\in\CO$,
because $R\Delta_\sigma=\Delta_\sigma(\Delta_{\sigma^{-1}} R\Delta_\sigma)$
and $\Delta_{\sigma^{-1}} R\Delta_\sigma\in\CO$.

Instead of $\CO$, we can use the basis $\{f_R\mid R\in\CO \}$ of~$\CP_1$,
consisting of the idempotents of~$\CP_1$ defined in Theorem~\ref{semi-simple}.
The group $\Sigma$ acts by conjugation on the set $\CO$ of all orders,
hence also on the set $\{f_R\mid R\in\CO \}$. We first record the following easy observation.

\result{Lemma} \label{conj-fR}
Let $R$ be an order and let $f_R$ be the corresponding idempotent of~$\CP_1$.
For every $\sigma\in\Sigma$,
$$\Delta_\sigma f_R\Delta_{\sigma^{-1}}= f_{\ls \sigma R} \,,$$
where $\ls\sigma R:=\Delta_\sigma R \Delta_{\sigma^{-1}}$.
\fresult

\pf
This follows immediately from the definition of~$f_R$ in Section~6.
\endpf

Since $\CP=\bigoplus_{\sigma\in\Sigma}\CP_\sigma$, this has a $k$-basis
$\{\Delta_\sigma f_R\mid \sigma\in\Sigma\,,\,R\in\CO\}$.
We now describe the product in $\CP$ with respect to this basis.

\result{Lemma} \label{product}
The product of basis elements of $\CP$ is given by~:
$$(\Delta_\tau f_S)(\Delta_\sigma f_R)=
\left\{\begin{array}{ll}
0    \;\;   &\text{if } S\neq\ls\sigma R \,, \\
\Delta_{\tau\sigma} f_R  \;\; &\text{if }\; S=\ls \sigma R \,,
\end{array}\right.$$
for all $S,R\in\CO$ and all $\tau,\sigma\in\Sigma$.
\fresult

\pf
$(\Delta_\tau f_S) (\Delta_\sigma f_R)=\Delta_\tau f_S  f_{\ls\sigma R} \Delta_\sigma$.
This is zero if $S\neq\ls\sigma R$. Otherwise we obtain
$\Delta_\tau f_{\ls\sigma R} \Delta_\sigma=\Delta_\tau \Delta_\sigma f_R
=\Delta_{\tau\sigma} f_R$.
\endpf

\result{Corollary} \label{proj}
Let $R$ be an order and let $f_R$ be the corresponding idempotent of~$\CP_1$.
The left ideal $\CP f_R$ has a $k$-basis
$\{\Delta_\sigma f_R\mid \sigma\in\Sigma\}$.
\fresult

\pf
We have $(\Delta_\tau f_S) f_R=\left\{\begin{array}{ll}
0    \;\;   &\text{if } S\neq\ls\sigma R \,, \\
\Delta_{\tau} f_R  \;\; &\text{if }\; S=\ls \sigma R \,,
\end{array}\right.$
and we know that the set $\{\Delta_\sigma f_R\mid \sigma\in\Sigma\}$ is part of the basis of~$\CP$.
\endpf

Now we introduce central idempotents in~$\CP$.
Let $R$ be an order, let $\Sigma_R$ be the stabilizer of~$R$ in~$\Sigma$,
and denote by $[\Sigma/\Sigma_R]$ a set of coset representatives.
By Lemma~\ref{conj-fR}, $\Sigma_R$ is also the stabilizer of~$f_R$ and we define
$$e_R=\sum_{\sigma\in[\Sigma/\Sigma_R]} \Delta_\sigma f_R\Delta_{\sigma^{-1}}
=\sum_{\sigma\in[\Sigma/\Sigma_R]} f_{\ls \sigma R} \,,$$
the sum of the $\Sigma$-orbit of~$f_R$.

\result{Lemma} \label{central-idempotents}
Let $[\Sigma\backslash\CO]$ be a set of representatives of the $\Sigma$-orbits in~$\CO$.
The set $\{ e_R \mid R\in [\Sigma\backslash\CO] \}$ is a set of orthogonal central idempotents of~$\CP$,
whose sum is $1_{\CP}=\Delta$.
\fresult

\pf
We compute 
$$\Delta_\tau f_S \,e_R= \!\!\!
\sum_{\sigma\in[\Sigma/\Sigma_R]}  \!\!\! \Delta_\tau f_S \,f_{\ls \sigma R}
=\left\{\begin{array}{ll}
0    \,  &\text{if $S$ does not belong to the orbit of~$R$} \,, \\
\Delta_\tau f_S  \, &\text{if }\; S=\ls \sigma R \,.
\end{array}\right.$$
On the other hand
$$e_R \,\Delta_\tau f_S=
\sum_{\sigma\in[\Sigma/\Sigma_R]}  f_{\ls \sigma R} \,\Delta_\tau f_S
=\sum_{\sigma\in[\Sigma/\Sigma_R]}  f_{\ls \sigma R} \,f_{\ls\tau S} \,\Delta_\tau \,.
$$
This is zero if $\ls\tau S$ does not belong to the $\Sigma$-orbit of~$\ls\sigma R$,
that is, if $S$ does not belong to the $\Sigma$-orbit of~$R$, while if $\ls\tau S=\ls\sigma R$,
then we get $f_{\ls\tau S} \,\Delta_\tau=\Delta_\tau f_S$. This shows that $e_R$ is central.

We know that $\{f_R \mid R\in\CO \}$ is a set of orthogonal idempotents with sum~1.
Since we have just grouped together the $\Sigma$-orbits,
it is clear that the set $\{ e_R \mid R\in [\Sigma\backslash\CO] \}$
is also a set of orthogonal idempotents of~$\CP$, whose sum is $1_{\CP}=\Delta$.
\endpf

It follows from Lemma~\ref{central-idempotents} that
$\CP\cong \prod_{R\in [\Sigma\backslash\CO]} \CP e_R$ and we have to understand the structure of each term.

\result{Theorem} \label{structure}
Let $R$ be an order on~$X$ and let $\Sigma_R$ be its stabilizer in~$\Sigma$. Then
$$\CP e_R\cong M_{|\Sigma:\Sigma_R|}(k\Sigma_R) \,,$$
a matrix algebra of size $|\Sigma:\Sigma_R|$ on the group algebra $k\Sigma_R$.
In other words
$$\CP\cong \prod_{R\in [\Sigma\backslash\CO]} M_{|\Sigma:\Sigma_R|}(k\Sigma_R) \,.$$
\fresult

\pf
By Corollary~\ref{proj}, the left ideal $\CP f_R$ is a free $k$-submodule of~$\CP$
spanned by the set $\{\Delta_\sigma f_R \mid \sigma\in\Sigma \}$.
The group $\Sigma_R$ acts on the right on this set,
because $f_R\Delta_h =\Delta_h f_R$ for every $h\in\Sigma_R$.
It follows that $\CP f_R$ is a free right $k\Sigma_R$-module with basis
$\{\Delta_\sigma f_R \mid \sigma\in[\Sigma/\Sigma_R] \}$.

Clearly, the left action of $\CP$ commutes with the right action of~$k\Sigma_R$.
The left action of $\CP$ on this free right $k\Sigma_R$-module induces a $k$-algebra map
$$\phi_R:\CP\longrightarrow M_{|\Sigma:\Sigma_R|}(k\Sigma_R) \,.$$
By Lemma~\ref{product}, $e_R$ acts as the identity on~$\CP f_R$,
while $e_S$ acts by zero if $S$ does not belong to the $\Sigma$-orbit of~$R$.
Therefore we get a $k$-algebra map
$$\phi_R:\CP e_R\longrightarrow M_{|\Sigma:\Sigma_R|}(k\Sigma_R) \,,$$
because $\phi_R(e_S)=0$ whenever $S$ does not belong to the $\Sigma$-orbit of~$R$.
Putting all these maps together, we obtain a $k$-algebra map
$$\phi \,=\prod_{R\in [\Sigma\backslash\CO]} \phi_R \;: \;\;
\CP\longrightarrow \prod_{R\in [\Sigma\backslash\CO]}
M_{|\Sigma:\Sigma_R|}(k\Sigma_R) \,.$$
Lemma~\ref{product} shows that, if $\rho,\sigma,\tau\in [\Sigma/\Sigma_R]$
and $g\in\Sigma_R$, we have
$$(\Delta_\tau f_R\Delta_g \Delta_{\rho^{-1}}) \cdot \Delta_\sigma f_R
=\left\{\begin{array}{ll}
0    \;\;   &\text{if } \rho\neq\sigma \,, \\
\Delta_\tau f_R \Delta_g \;\; &\text{if }\; \rho=\sigma \,,
\end{array}\right.$$
using the fact that $f_R \Delta_g= \Delta_g f_R$.
This means that $\phi_R(\Delta_\tau f_R\Delta_g\Delta_{\rho^{-1}})$ is the elementary matrix
with a single nonzero entry equal to~$\Delta_g$ in position $(\tau,\rho)$.
Moreover, we also have $\phi_S(\Delta_\tau f_R\Delta_g\Delta_{\rho^{-1}})=0$
whenever $S$ does not belong to the $\Sigma$-orbit of~$R$.
Therefore the map $\phi$ is surjective.

Finally, we prove that $\phi$ is an isomorphism.
It suffices to do this in the case where $k=\Z$, because all the algebras are defined over~$\Z$
(that is, they are obtained by extending scalars from $\Z$ to~$k$)
and the algebra map~$\phi$ is also defined over~$\Z$.
Now if $k=\Z$, then all algebras under consideration are finitely generated free $\Z$-modules
and we know that the map $\phi$ is surjective.
So it suffices to show that the source and the target of~$\phi$ have the same rank as $\Z$-modules.
The rank of~$\CP$ is $|\Sigma|\Card(\CO)$. On the other hand,
$${\rm rank}(M_{|\Sigma:\Sigma_R|}(k\Sigma_R))
=|\Sigma:\Sigma_R|^2|\Sigma_R| =|\Sigma:\Sigma_R| \, |\Sigma| \,.$$
Summing over $R\in [\Sigma\backslash\CO]$, we obtain
$$\sum_{R\in [\Sigma\backslash\CO]} |\Sigma:\Sigma_R| \, |\Sigma|
=|\Sigma| \sum_{R\in [\Sigma\backslash\CO]} \Card(\text{orbit of }R) =|\Sigma|\Card(\CO) \,,$$
as was to be shown.
\endpf

\result{Remark}\label{Morita} {\rm Since a matrix algebra $M_r(A)$ is Morita equivalent to~$A$
(for any $k$-algebra~$A$), it follows from Theorem~\ref{structure} that
the algebra $\CP$ is Morita equivalent to a product of group algebras,
namely $B=\prod_{R\in[\Sigma\backslash\CO]} k\Sigma_R$.
The bimodule which provides the Morita equivalence is $M=\bigoplus_{R\in[\Sigma\backslash\CO]} \CP f_R$,
which is clearly a left $\CP$-module by left multiplication,
and a right module for each group algebra $k\Sigma_R$,
acting by right multiplication on the summand $\CP f_R$, and
acting by zero on the other summands $\CP f_S$, where $S\neq R$ in $[\Sigma\backslash\CO]$.
Notice that $\CP f_R$ is the bimodule appearing in the proof of Theorem~\ref{structure}.

The bimodule inducing the inverse Morita equivalence is
$M^\vee=\bigoplus_{R\in[\Sigma\backslash\CO]}  f_R\CP$. Indeed we, obtain
$$M\otimes_B M^\vee
\cong \bigoplus_{S\in[\Sigma\backslash\CO]}\bigoplus_{R\in[\Sigma\backslash\CO]} \CP f_Sf_R \CP
=\bigoplus_{R\in[\Sigma\backslash\CO]} \CP f_R \CP=\CP$$
and on the other hand
$$\begin{array}{rl}
M^\vee\otimes_\CP M \cong& \displaystyle \bigoplus_{S\in[\Sigma\backslash\CO]}
\bigoplus_{R\in[\Sigma\backslash\CO]} f_S\CP f_R
=\bigoplus_{R\in[\Sigma\backslash\CO]} f_R\CP f_R \\
=&\displaystyle\bigoplus_{R\in[\Sigma\backslash\CO]}k\Sigma_R{\cdot}f_R
\cong \bigoplus_{R\in[\Sigma\backslash\CO]} k\Sigma_R = B\,.
\end{array}$$
Each module $\CP f_R$ has rank $|\Sigma|$, because it has a $k$-basis
$\{\Delta_\sigma f_R \mid \sigma\in\Sigma\}$,
but we view it as a free right $k\Sigma_R$-module of rank~$|\Sigma:\Sigma_R|$.
}
\fresult


\vspace{-.1cm}
\section{Simple modules for the essential algebra}

\noindent
By standard commutative algebra, any simple $\CE$-module is actually a module over $k/m\otimes_k\CE$,
where $m$ is a maximal ideal of~$k$.
Replacing $k$ by the quotient $k/m$, we assume from now on that $k$ is a field.
Let $\CE$ be the essential algebra of Section~5 and let $\CP=\CE/N$ be the algebra of permuted orders.
Since $N$ is a nilpotent ideal and since nilpotent ideals act by zero on simple modules,
any simple $\CE$-module can be viewed as a simple $\CP$-module.
So we work with~$\CP$ and we wish to describe all simple left $\CP$-modules.

There is a general procedure for constructing all simple modules for the algebra of a semigroup~$S$,
using equivalence classes of maximal subgroups of~$S$,
see Theorem~5.33 in~\cite{CP}, or Section~3 of \cite{HK} for a short presentation.
But our previous results allow for a very direct and easy approach,
so we do not need to follow the method of~\cite{CP}.

First notice that the simple $\CP_1$-modules are easy to describe,
because $\CP_1$ is a product of copies of~$k$ (by Theorem~\ref{semi-simple}).
More precisely, $\CP_1\cong \prod_{R\in\CO} k\cdot f_R$
and each one-dimensional space $k\cdot f_R$ is a simple $\CP_1$-module
(where $R$ runs through the set~$\CO$ of all orders).

\result{Theorem} \label{simple}
Assume that $k$ is a field.
Let $\CW$ be the set of all pairs $(R,V)$, where $R$ is an order on~$X$
and $V$ is a simple $k\Sigma_R$-module up to isomorphism.
The group $\Sigma$ acts on~$\CW$ via $\ls\sigma(R,V):=(\ls\sigma R,\ls\sigma V)$,
where $\ls\sigma R=\Delta_\sigma R\Delta_{\sigma^{-1}}$ and $\ls\sigma V$ is the conjugate module,
a module for the group algebra $k\Sigma_{\ls\sigma R}=k[\sigma \Sigma_R\,\sigma^{-1}]$.
\begin{itemize}
\item[(a)] The set of isomorphism classes of simple $\CP$-modules is parametrized by
the set $\Sigma\backslash\CW$ of $\Sigma$-conjugacy classes of pairs $(R,V)\in\CW$.
\item[(b)] The simple module corresponding to~$(R,V)$ under the parametrization of part~(a) is
$$S_{R,V} = W_R\otimes_k V \,,$$
where $W_R$ is the unique (up to isomorphism) simple module for the matrix algebra
$M_{|\Sigma:\Sigma_R|}(k)$ and $W_R\otimes V$ is viewed as a module for the algebra
$$M_{|\Sigma:\Sigma_R|}(k) \otimes_k k\Sigma_R \cong 
M_{|\Sigma:\Sigma_R|}(k\Sigma_R) \,,$$
which is one of the factors of the decomposition of~$\CP$ in Theorem~\ref{structure}.
\item[(c)] The simple $\CP$-module $S_{R,V}$ is also isomorphic to $\CP f_R\otimes_{k\Sigma_R}V$,
with its natural structure of $\CP$-module under left multiplication.
\item[(d)] The simple $\CP$-module $S_{R,V}$ has dimension $|\Sigma:\Sigma_R|\cdot\dim(V)$.
\end{itemize}
\fresult

\pf
By Theorem~\ref{structure}, any simple $\CP$-module is a simple module for one of the factors
$M_{|\Sigma:\Sigma_R|}(k\Sigma_R)$, where $R$ belongs to
a set $[\Sigma\backslash\CO]$ of representatives of $\Sigma$-orbits in~$\CO$.
In view of the isomorphism
$$M_{|\Sigma:\Sigma_R|}(k\Sigma_R) \cong 
M_{|\Sigma:\Sigma_R|}(k) \otimes_k k\Sigma_R \,,$$
any such simple module is isomorphic to a tensor product $W_R\otimes_k V$ as in the statement.
This proves (a) and (b).

By Theorem~\ref{structure}, $\CP e_R\cong M_{|\Sigma:\Sigma_R|}(k\Sigma_R)$
and its identity element $e_R$ decomposes as a sum of orthogonal idempotents
$e_R=\sum_{\sigma\in[\Sigma/\Sigma_R]} f_{\,^\sigma \! R}\,$.
Cutting by the idempotent $f_R$, we obtain the left ideal $\CP f_R$,
which is a free right $k\Sigma_R$-module,
isomorphic to the space of column vectors with coefficients in~$k\Sigma_R$.
Now $\CP f_R$ is the bimodule providing the Morita equivalence between
$M_{|\Sigma:\Sigma_R|}(k\Sigma_R)$ and $k\Sigma_R$ (see Remark~\ref{Morita}).
Therefore, for any simple left $k\Sigma_R$-module~$V$, the corresponding simple module for
$\CP e_R\cong M_{|\Sigma:\Sigma_R|}(k\Sigma_R)$ is the left $\CP$-module
$\CP f_R\otimes_{k\Sigma_R}V$.
Since $\CP f_R$ is the space of column vectors with coefficients in~$k\Sigma_R$, 
while $W_R$ is the space of column vectors with coefficients in~$k$, we get
$\CP f_R\cong W_R\otimes_k k\Sigma_R$. Therefore our simple $\CP$-module~is
$$\CP f_R\otimes_{k\Sigma_R}V\cong W_R\otimes_k k\Sigma_R\otimes_{k\Sigma_R}V
\cong W_R\otimes_k V \cong S_{R,V} \,,$$
proving (c).

Finally, the dimension is
$$\dim(S_{R,V})=\dim(W_R\otimes_k V)=\dim(W_R)\cdot \dim(V)=|\Sigma:\Sigma_R|\cdot\dim(V) \,,$$
proving~(d).
\endpf

\result{Example} \label{Delta}
{\rm Consider the trivial order~$\Delta$. Then $\Sigma_\Delta=\Sigma$ and the matrix algebra reduces to
$$M_{|\Sigma:\Sigma_\Delta |}(k\Sigma_\Delta) \cong M_1(k)\otimes_k k\Sigma \cong k\Sigma \,.$$
The simple module $W_\Delta$ for the algebra $M_1(k)$ is just~$k$ and
the simple module $S_{\Delta,V}=W_\Delta\otimes_k V\cong V$ is just a simple $k\Sigma$-module.
In that case, the central idempotent~$e_S$ of~$\CP$ acts by zero on~$V$ for any order $S\neq\Delta$,
hence $f_S$ too (because $f_S\,e_S=f_S$).
Then $R=\sum_{R\subseteq S} f_S$ also acts by zero for any order $R\neq\Delta$.
For any essential reflexive relation $Q$ with $Q\neq\Delta$,
the action of $Q$ is equal to the action of~$\overline Q$
(because $Q-\overline Q$ belongs to the nilpotent ideal~$N$),
and therefore $Q$ also acts by zero on~$V$. Then so does the action
of the essential relation $\Delta_\sigma Q$ containing the permutation~$\Delta_\sigma$.
This shows that the simple modules $S_{\Delta,V}\cong V$ are just
the modules for $k\Sigma$ viewed as a quotient algebra as in Lemma~\ref{group-algebra}.
}
\fresult

\result{Example} \label{total-orders}
{\rm Consider a total order~$T$. Then $\Sigma_T=\{\Id\}$ and the matrix algebra reduces to
$$M_{|\Sigma:\Sigma_T |}(k\Sigma_T) \cong M_{n!}(k)\otimes_k k \cong M_{n!}(k) \,.$$
In that case, there is unique simple $k\Sigma_T$-module, namely $V=k$, the trivial module for the trivial group.
We obtain the single simple module $S_{T,k}=W_T\otimes_k k\cong W_T$ for the algebra $M_{n!}(k)$.
Equivalently, with the approach of part~(c) of Theorem~\ref{simple},
we have $f_T=T$ (by maximality of~$T$ in~$\CO$) and so
$$S_{T,k}=\CP f_T\otimes_{k\Sigma_T}V=\CP T\otimes_{k}k\cong \CP T \,.$$
So we obtain just the left ideal $\CP f_T=\CP T$, which turns out to be simple in that case.
But it is also the left ideal~$L$ appearing in Proposition~\ref{n!}.
So we have recovered the simple module of Proposition~\ref{n!}.
}
\fresult

We also mention another byproduct of Theorem~\ref{structure}.

\result{Theorem} \label{Maschke}
If the characteristic of the field $k$ is zero or $>n$, then $\CP$ is a semi-simple $k$-algebra.
\fresult

\pf
It suffices to see that each factor in the decomposition of Theorem~\ref{structure} is semi-simple.
Now we have the isomorphism
$$M_{|\Sigma:\Sigma_R|}(k\Sigma_R) \cong 
M_{|\Sigma:\Sigma_R|}(k) \otimes_k k\Sigma_R \,,$$
and $M_{|\Sigma:\Sigma_R|}(k)$ is a simple algebra.
Moreover the group algebra $k\Sigma_R$ is semi-simple by Maschke's theorem,
because the characteristic of~$k$ does not divide the order of the group $\Sigma_R$,
by assumption. The result follows.
\endpf

Every simple $\CP$-module $S_{R,V}$ is a simple $\CR$-module, because of the successive quotients
$\CR\to\CE\to\CP$.
We now give a direct description of the action on~$S_{R,V}$ of an arbitrary relation in~$\CR$.
Since $S_{R,V}\cong\CP f_R\otimes_{k\Sigma_R}V$, it suffices to describe the action on~$\CP f_R$,
and for this we can work again with an arbitrary commutative base ring~$k$.
Recall that $\CP f_R$ has a basis $\{\Delta_\sigma f_R \,|\, \sigma\in\Sigma \}$.

\result{Proposition} 
Let $k$ be a commutative ring.
Let $R$ be an order on~$X$ and let $Q$ be an arbitrary relation (in the $k$-algebra~$\CR$).
The action of $Q$ on $\CP f_R$ is described on the basis elements as follows~:
$$Q\cdot \Delta_\sigma f_R =
\left\{\begin{array}{cl}
\Delta_{\tau\sigma}f_R &\text{ if }\,\exists\, \tau\in\Sigma \text{ such that }
\Delta\subseteq \Delta_{\tau^{-1}}Q \subseteq \ls\sigma R \,, \\
0 & \text{ otherwise }.
\end{array}\right.$$
\fresult

\pf
Suppose first that $S$ is an order. By Lemma~\ref{product}, the action of $f_S$ is given by
$$f_S\cdot \Delta_\sigma f_R =
\left\{\begin{array}{cl}
\Delta_{\sigma}f_R &\text{ if }\,S = \ls\sigma R \,, \\
0 & \text{ otherwise }.
\end{array}\right.$$
Now $\displaystyle S=\sum_{\scriptstyle T\in\CO \atop\scriptstyle S\subseteq T}  f_T$
and the action of $f_T$ is nonzero only if $T=\ls\sigma R$.
So we obtain the action of~$S$ as follows~:
$$S\cdot \Delta_\sigma f_R =
\left\{\begin{array}{cl}
\Delta_{\sigma}f_R &\text{ if }\,S \subseteq \ls\sigma R \,, \\
0 & \text{ otherwise }.
\end{array}\right.$$
Next we suppose that $S$ is reflexive and that its transitive closure $\overline S$ is an order.
Since $S-\overline S$ belongs to the nilpotent ideal $N$ of Section~5,
which acts by zero because $\CP=\CE/N$, the action of $S$ coincides with the action of~$\overline S$.
Moreover, the condition $\overline S \subseteq \ls\sigma R$ is equivalent to $S \subseteq \ls\sigma R$,
because $\ls\sigma R$ is transitive. Therefore the action of~$S$ is the following~:
$$S\cdot \Delta_\sigma f_R =
\left\{\begin{array}{cl}
\Delta_{\sigma}f_R &\text{ if }\,S \subseteq \ls\sigma R \,, \\
0 & \text{ otherwise }.
\end{array}\right.$$
Now suppose that $S$ is reflexive and that $\overline S$ is not an order.
Then $\overline S$ is inessential by Proposition~\ref{preorder}, hence zero in~$\CE$.
So $\overline S$ acts by zero, and since $S-\overline S$ acts by zero, the action of $S$ is also zero.
On the other hand, $S$ cannot be contained in $\ls\sigma R$,
otherwise $\overline S\subseteq \ls\sigma R$, which would force $\overline S$ to be an order
since $\ls\sigma R$ is an order.
Therefore the condition $S \subseteq \ls\sigma R$ is never satisfied in that case.
So the previous formula still holds, because we have zero on both sides~:
$$S\cdot \Delta_\sigma f_R =
\left\{\begin{array}{cl}
\Delta_{\sigma}f_R &\text{ if }\,S \subseteq \ls\sigma R \,, \\
0 & \text{ otherwise }.
\end{array}\right.$$
Now suppose that $Q$ contains a permutation $\Delta_\tau$. Then $S=\Delta_{\tau^{-1}} Q$ is reflexive and
$Q=\Delta_\tau S$. Thus the action of $Q$ is~:
$$\begin{array}{rl}
Q\cdot \Delta_\sigma f_R &=
\left\{\begin{array}{cl}
\Delta_\tau\Delta_{\sigma}f_R &\text{ if }\,S \subseteq \ls\sigma R \,, \\
0 & \text{ otherwise },
\end{array}\right. \\
&=\left\{\begin{array}{cl}
\Delta_{\tau\sigma}f_R &\text{ if }\,\Delta_{\tau^{-1}} Q \subseteq \ls\sigma R \,, \\
0 & \text{ otherwise },
\end{array}\right. \\
&=\left\{\begin{array}{cl}
\Delta_{\tau\sigma}f_R &\text{ if }\,\Delta \subseteq\Delta_{\tau^{-1}} Q \subseteq \ls\sigma R \,, \\
0 & \text{ otherwise }.
\end{array}\right.
\end{array}$$
The last equality holds because the condition $\Delta_{\tau^{-1}}Q \subseteq \ls\sigma R$ is equivalent to
$\Delta \subseteq\Delta_{\tau^{-1}}Q \subseteq \ls\sigma R$, since $S=\Delta_{\tau^{-1}} Q$ is reflexive.
This proves the result for such a relation~$Q$.

Finally if $Q$ does not contain a permutation, then $Q$ is inessential by Theorem~\ref{contain},
hence acts by zero. On the other hand the condition that there exists $\tau\in\Sigma$ such that
$\Delta \subseteq\Delta_{\tau^{-1}}Q$ cannot be satisfied since $Q$ does not contain a permutation.
Therefore the previous formula still holds, because we have zero on both sides~:
$$Q\cdot \Delta_\sigma f_R =
\left\{\begin{array}{cl}
\Delta_{\tau\sigma}f_R &\text{ if }\,\Delta \subseteq\Delta_{\tau^{-1}} Q \subseteq \ls\sigma R \,, \\
0 & \text{ otherwise }.
\end{array}\right.$$
This proves the result in all cases.
\endpf

\result{Remark} \label{finite-group}
{\rm In the description of the algebra~$\CP$ (Theorem~\ref{structure}) and in the description of
its simple modules (Theorem~\ref{simple}), we may wonder which groups
appear as $\Sigma_R$ for some order~$R$. The answer is that the group $\Sigma_R$ is arbitrary.
More precisely, by a theorem of Birkhoff~\cite{Bi} (and further improvements by Thornton~\cite{Tho}
and Barmak--Minian~\cite{BaM}), any finite group is isomorphic to $\Sigma_R$
for some order~$R$ on a suitable finite set~$X$.
However, for a given finite set~$X$,
it is not clear which isomorphism classes of groups $\Sigma_R$ appear.

Another question is to determine whether or not the simple modules of Theorem~\ref{simple}
are absolutely simple. But again this depends on the group $\Sigma_R$,
because the field $k$ may or may not be a splitting field for the group algebra~$k\Sigma_R$.}
\fresult


\vspace{-.1cm}
\section{A branching rule}

\noindent
In this section, we let $X=\{1,\ldots,n\}$ for simplicity.
In order to let $n$ vary, we use a superscript $(n)$ for all objects depending on~$n$,
such as $X^{(n)}$ for the set~$X$,
$\Sigma^{(n)}$ for the symmetric group on~$X^{(n)}$,
$\CO^{(n)}$ for the set of all orders on~$X^{(n)}$,
$\CP_1^{(n)}$ for the algebra of orders on~$X^{(n)}$,
$\CP^{(n)}$ for the algebra of permuted orders on~$X^{(n)}$, etc.

In the representation theory of the symmetric group $\Sigma^{(n)}$,
there are well-known branching rules,
describing the restriction of simple modules to the subgroup~$\Sigma^{(n-1)}$
of all permutations of $X^{(n-1)}$, and the induction of simple modules from $\Sigma^{(n-1)}$
to~$\Sigma^{(n)}$. In a similar fashion,
working again over an arbitrary commutative base ring~$k$,
we will describe how modules for $\CP^{(n-1)}$ behave under induction to~$\CP^{(n)}$.
For this we need to view the former as a subalgebra of the latter. We first define
$$\phi:\CP_1^{(n-1)}\longrightarrow \CP_1^{(n)} \,,\qquad \phi(R)=R\cup\{ (n,n)\} \,,$$
for any order~$R$ on $X^{(n-1)}$. It is clear that $\phi(R)$ is an order on~$X^{(n)}$.
Since $\Sigma^{(n-1)}$ is a subgroup of~$\Sigma^{(n)}$ (by fixing the last letter~$n$),
the map $\phi$ clearly extends to an injective algebra homomorphism
$\phi:\CP^{(n-1)}\longrightarrow \CP^{(n)}$.

Now we want to compute the image under~$\phi$ of the idempotents $f_R^{(n-1)}$.
For a given order $R$ on~$X^{(n-1)}$, we define
$$\CS_R= \{ S\in \CO^{(n)} \mid S \cap (X^{(n-1)}\times X^{(n-1)}) =R\} \,.$$

\result{Lemma} \label{n-1,n}
Let $R$ be an order on~$X^{(n-1)}$
and let $f_R^{(n-1)}$ be the corresponding idempotent of~$\CP_1^{(n-1)}$.

\begin{itemize}
\item[(a)] If $S$ is an order on~$X^{(n)}$
and if $f_S^{(n)}$ is the corresponding idempotent of~$\CP_1^{(n)}$, then
$$\phi(f_R^{(n-1)})\cdot f_S^{(n)} =\left\{\begin{array}{cl}
f_S^{(n)} &\text{ if }\,S \in\CS_R  \,, \\
0 & \text{ otherwise }.
\end{array}\right.$$

\item[(b)] $\;\phi(f_R^{(n-1)})\,= \displaystyle\sum_{S \in\CS_R} f_S^{(n)} $.

\end{itemize}

\fresult

\pf
We have $f_R^{(n-1)}
=\displaystyle \sum_{\scriptstyle Y\in\CO^{(n-1)} \atop\scriptstyle R\subseteq Y} \mu(R,Y) Y \;$
and $\;\phi(Y)=Y\cup\{ (n,n)\}$. Therefore, using part~(a) of Theorem~\ref{inversion}, we obtain
$$\begin{array}{cl}
\phi(f_R^{(n-1)})\cdot f_S^{(n)}
&=\displaystyle\sum_{\scriptstyle Y\in\CO^{(n-1)} \atop\scriptstyle R\subseteq Y}
\mu(R,Y) (Y\cup\{ (n,n)\}) \cdot f_S^{(n)} \\
&=\displaystyle\Big(\sum_{\scriptstyle Y\in\CO^{(n-1)} \atop\scriptstyle R\subseteq Y, \,Y\cup\{ (n,n)\} \subseteq S}
\mu(R,Y) \Big) f_S^{(n)} \\
&=\displaystyle\Big(\sum_{\scriptstyle Y\in\CO^{(n-1)} \atop\scriptstyle
R\subseteq Y \subseteq S \cap (X^{(n-1)}\times X^{(n-1)})} \mu(R,Y) \Big) f_S^{(n)} \,.
\end{array}$$
We get zero if $R\not\subseteq S \cap (X^{(n-1)}\times X^{(n-1)})$
and also if $R\subset S \cap (X^{(n-1)}\times X^{(n-1)})$ (by the definition of the M\"obius function).
If now $R= S \cap (X^{(n-1)}\times X^{(n-1)})$, that is, if $S \in\CS_R$,
then the sum reduces to $\mu(R,R)=1$ and we obtain $f_S^{(n)}$, proving~(a).

Now we have
$$\phi(f_R^{(n-1)})=\sum_{S\in \CO^{(n)}}\phi(f_R^{(n-1)})\cdot f_S^{(n)}
=\sum_{S\in \CS_R} f_S^{(n)} \,,$$
proving (b).
\endpf

\result{Theorem} \label{branching}
Let $R$ be an order on~$X^{(n-1)}$,
let $f_R^{(n-1)}$ be the corresponding idempotent of~$\CP_1^{(n-1)}$,
and let $V$ be a $k\Sigma_R^{(n-1)}$-module.
Then, inducing to~$\CP^{(n)}$ the $\CP^{(n-1)}$-module
$S_{R,V}=\CP^{(n-1)}f_R^{(n-1)}\otimes_{k\Sigma_R^{(n-1)}}V$, we obtain
$$\begin{array}{rl}
\CP^{(n)}&\!\!\!\otimes_{\CP^{(n-1)}} \big(\CP^{(n-1)}f_R^{(n-1)}\otimes_{k\Sigma_R^{(n-1)}}V\big) \\
&\displaystyle \cong \bigoplus_{S\in \CS_R}
\CP^{(n)}f_S^{(n)} \otimes_{k\Sigma_S^{(n)}}
\Ind_{\Sigma_R^{(n-1)}\cap\Sigma_S^{(n)}}^{\Sigma_S^{(n)}}
\Res_{\Sigma_R^{(n-1)}\cap\Sigma_S^{(n)}}^{\Sigma_R^{(n-1)}} V \,.
\end{array}$$
\fresult

\pf Using Lemma~\ref{n-1,n}, we obtain
$$\begin{array}{rl}
\CP^{(n)}&\!\!\! \otimes_{\CP^{(n-1)}} \big(\CP^{(n-1)}f_R^{(n-1)}\otimes_{k\Sigma_R^{(n-1)}}V\big) \\
&=\CP^{(n)}\phi(f_R^{(n-1)})\otimes_{k\Sigma_R^{(n-1)}}V \\
&\cong\displaystyle  \bigoplus_{S\in\CS_R} \CP^{(n)}f_S^{(n)} \otimes_{k\Sigma_R^{(n-1)}}V \\
&\displaystyle \cong \bigoplus_{S\in \CS_R}
\CP^{(n)}f_S^{(n)} \otimes_{k\Sigma_S^{(n)}} k\Sigma_S^{(n)}
\otimes_{k[\Sigma_R^{(n-1)}\cap\Sigma_S^{(n)}]} k[\Sigma_R^{(n-1)}\cap\Sigma_S^{(n)}]
\otimes_{k\Sigma_R^{(n-1)}}V \\
&\displaystyle \cong \bigoplus_{S\in \CS_R}
\CP^{(n)}f_S^{(n)} \otimes_{k\Sigma_S^{(n)}}
\Ind_{\Sigma_R^{(n-1)}\cap\Sigma_S^{(n)}}^{\Sigma_S^{(n)}}
\Res_{\Sigma_R^{(n-1)}\cap\Sigma_S^{(n)}}^{\Sigma_R^{(n-1)}} V \,,
\end{array}$$
proving the result.
\endpf

Assume for simplicity that the base ring $k$ is a field of characteristic zero and let
$V$ be a simple $k\Sigma_R^{(n-1)}$-module.
The $k\Sigma_S^{(n)}$-module
$\Ind_{\Sigma_R^{(n-1)}\cap\Sigma_S^{(n)}}^{\Sigma_S^{(n)}}
\Res_{\Sigma_R^{(n-1)}\cap\Sigma_S^{(n)}}^{\Sigma_R^{(n-1)}} V$ is a direct sum of simple modules~$W$,
and each $W$ gives rise to a simple $\CP^{(n)}$-module
$\CP^{(n)}f_S^{(n)} \otimes_{k\Sigma_S^{(n)}} W$.
Moreover, every such simple $\CP^{(n)}$-module occurs with multiplicities,
appearing for instance whenever we have
$\sigma$ running in a set of representatives of cosets
$[\Sigma_R^{(n-1)}/ \Sigma_R^{(n-1)}\cap\Sigma_S^{(n)}]$.
For any such $\sigma$, we have $\ls \sigma S\in\CS_R$ and also $\ls \sigma V\cong V$,
because $V$ is a $k\Sigma_R^{(n-1)}$-module and $\sigma\in\Sigma_R^{(n-1)}$.
Therefore the corresponding term in the direct sum is
$$\CP^{(n)}f_{\ls \sigma S}^{(n)} \otimes_{k\Sigma_R^{(n-1)}}V
\cong \CP^{(n)}f_{\ls \sigma S}^{(n)} \otimes_{k\Sigma_R^{(n-1)}}\ls\sigma V \,,$$
but this gives rise to the same simple $\CP^{(n)}$-modules as the ones coming from $S$,
by Theorem~\ref{simple}. Thus the multiplicity of these simple $\CP^{(n)}$-modules is at least
$|\Sigma_R^{(n-1)}/ \Sigma_R^{(n-1)}\cap\Sigma_S^{(n)}|$.


\vspace{-.1cm}
\section{Appendix on M\"obius inversion}

\noindent
In this appendix, we prove a general result on M\"obius inversion involving idempotents in a ring.
This was already used by the first author in other contexts (see Section~6.2 of~\cite{Bo})
and can be of independent interest.

Let $(P,\leq)$ be a finite lattice. Write $0$ for the minimal element of $P$
and write $x\vee y$ for the least upper bound of $x$ and $y$ in~$P$.

\result{Theorem} \label{inversion}
Let $P$ be a finite lattice.
Let $\{g_x \,|\, x\in P\}$ be a family of elements in a ring $A$ such that $g_0=1$ and
$g_x g_y=g_{x\vee y}$ for all $x,y\in P$.
For every $x\in P$, define
$$f_x= \sum_{\scriptstyle y\in P \atop\scriptstyle x\leq y} \mu(x,y) g_y \,,$$
where $\mu$ denotes the M\"obius function of the poset~$P$.

\begin{itemize}
\item[(a)] For all $x,y\in P$, we have $g_zf_x=f_x g_z= \left\{\begin{array}{ll}
f_x\; &\text{ if } z\leq x \,, \\
\,0 \; &\text{ if } z\not\leq x \,. \\
\end{array}\right.$

\item[(b)] The set $\{f_x \,|\, x\in P\}$ is a set of mutually orthogonal idempotents in~$P$ whose sum is~$1$.

\end{itemize}
\fresult

\medskip
\noindent
Note that our assumption implies that every $g_x$ is idempotent, because $x\vee x=x$.

\pf
By M\"obius inversion, we have
$$g_x=\sum_{\scriptstyle y\in P \atop\scriptstyle x\leq y} f_y \,,$$
and in particular $1=g_0=\sum_{\scriptstyle y\in P} f_y$.

Next we compute products. If $x,z\in P$, then
$$\begin{array}{rl}
f_x g_z =&\displaystyle
\Big( \sum_{\scriptstyle y\in P \atop\scriptstyle x\leq y} \mu(x,y)g_y \Big)  g_z 
=\displaystyle\sum_{\scriptstyle y\in P \atop\scriptstyle x\leq y} \mu(x,y)\; g_y  g_z
=\displaystyle\sum_{\scriptstyle y\in P \atop\scriptstyle x\leq y} \mu(x,y)\; g_{y\vee z} \\
=&\displaystyle\sum_{\scriptstyle w\in P \atop\scriptstyle x\vee z\leq w}
\Big(\sum_{\scriptstyle y\in P \atop\scriptstyle x\leq y,\;y \vee z=w} \mu(x,y)\Big) g_w \,.
\end{array}
$$
Note that $g_zf_x=f_x g_z$ because $g_zg_y=g_{z\vee y}=g_{y\vee z}=g_y  g_z$.
If $x$ is strictly smaller than~$x\vee z$, then the inner sum runs over the set of all elements~$y$
in the interval $[x,w]:=\{v\in P \mid x\leq v\leq w \}$ such that $y\vee(x\vee z)=w$. But we have
$$\sum_{\scriptstyle x\leq y \atop\scriptstyle y\vee z=w} \mu(x,y)
=\sum_{\scriptstyle x\leq y \atop\scriptstyle y\vee(x\vee z)=w} \mu(x,y) =0 \,,$$
by a well-known property of the M\"obius function (Corollary~3.9.3 in~\cite{St}).
Thus $f_x g_z=0$ if $x$ is strictly smaller than~$x\vee z$, that is, if $z\not\leq x$.

If now $x=x\vee z$, that is, $z\leq x$, we get $y=y\vee z$ (because $z\leq x\leq y$),
hence $y=w$, so that the inner sum has a single term for $y=w$. In that case, we get
$$f_x g_z=\sum_{\scriptstyle w\in P \atop\scriptstyle x\leq w} \mu(x,w)\,g_w = f_x \,.$$
Therefore
$$f_x g_z= \left\{\begin{array}{ll}
f_x\quad &\text{ if } z\leq x \,, \\
\,0 \quad &\text{ if } z\not\leq x \,, \\
\end{array}\right.$$
proving (a).

If now $x,u\in P$, then
$$f_x f_u
=\displaystyle \sum_{\scriptstyle y\in P \atop\scriptstyle u\leq y} \mu(u,y) \; f_x  g_y
=\displaystyle \sum_{\scriptstyle y\in P \atop\scriptstyle u\leq y\leq x}  \mu(u,y)\,f_x \,.
$$
If $u\not\leq x$, the sum is empty and we get zero. If $u< x$, then
$\displaystyle \sum_{u\leq y\leq x}  \mu(u,y)=0$
by the very definition of the M\"obius function.
This shows that $f_x f_u =0$ if $u\neq x$. Finally, if $u=x$, then we get
$f_x f_x = f_x$, thus $f_x$ is idempotent, and the proof is complete.
\endpf

\result{Corollary} \label{below}
Let $P$ be a finite lattice. Write $t$ for the maximal element of $P$
and write $x\wedge y$ for the greatest lower bound of $x$ and $y$ in~$P$.
Let $\{g_x \,|\, x\in P\}$ be a family of elements in a ring $A$ such that $g_t=1$ and
$g_x g_y=g_{x\wedge y}$ for all $x,y\in P$.
For every $x\in P$, define
$$f_x= \sum_{\scriptstyle y\in P \atop\scriptstyle y\leq x} \mu(y,x) g_y \,,$$
where $\mu$ denotes the M\"obius function of the poset~$P$.
Then the set $\{f_x \,|\, x\in P\}$ is a set of mutually orthogonal idempotents in~$P$ whose sum is~$1$.
\fresult

\pf
This follows from Theorem~\ref{inversion} by using the opposite ordering on~$P$.
\endpf


\bigskip
\noindent
Serge Bouc, CNRS-LAMFA, Universit\'e de Picardie - Jules Verne,\\
33, rue St Leu, F-80039 Amiens Cedex~1, France.\\
{\tt serge.bouc@u-picardie.fr}

\medskip
\noindent
Jacques Th\'evenaz, Section de math\'ematiques, EPFL, \\
Station~8, CH-1015 Lausanne, Switzerland.\\
{\tt Jacques.Thevenaz@epfl.ch}

\end{document}